\documentclass[thmsa,11pt]{article}
\usepackage{amsmath,amsthm,amssymb,amsfonts}
\usepackage{graphicx}
\usepackage[center]{caption}
\usepackage{subcaption}
\usepackage[noadjust]{cite}

\newtheorem{theorem}{Theorem}[section]

\newtheorem{lemma}[theorem]{Lemma}
\newtheorem{proposition}[theorem]{Proposition}

\newtheorem{remark}[theorem]{Remark}
\newtheorem{example}[theorem]{Example}

\begin{document}

\title{Multiple nonradial solutions for a nonlinear elliptic problem with singular
and decaying radial potential}
\author{{\small \smallskip }Sergio Rolando\smallskip \\
\textit{\small Dipartimento di Matematica e Applicazioni}\\
\textit{{\small Universit\`{a} degli Studi di Milano bicocca}}, 
\textit{\small Via Cozzi 53, 20125 Milano}, \textit{\small Italy}\\
\textit{\small e-mail: }{\small sergio.rolando@unito.it}}
\date{}
\maketitle

\begin{abstract}
Many existence and nonexistence results are known for nonnegative \emph{radial} 
solutions $u\in D^{1,2}(\mathbb{R}^{N})\cap L^{2}(\mathbb{R}^{N},\left|x\right| ^{-\alpha }dx)$ 
to the equation 
\[
-\triangle u+\displaystyle\frac{A}{\left| x\right| ^{\alpha }}u=f\left( u\right) \quad 
\textrm{in }\mathbb{R}^{N},\quad N\geq 3,\quad A,\alpha >0, 
\]
with nonlinearites satisfying $\left| f\left( u\right) \right| \leq \left(\mathrm{const.}\right) u^{p-1}$ 
for some $p>2$. Existence of \emph{nonradial} solutions, by contrast, is known only for $N\geq 4$, $\alpha =2$, 
$f\left( u\right) =u^{(N+2)/(N-2)}$ and $A$ large enough. Here we show that the equation has
multiple nonradial solutions as $A\rightarrow +\infty$ for $N\geq 4$, $2/(N-1)<\alpha <2N-2$, $\alpha\neq 2$,
and nonlinearities satisfying suitable assumptions. Our argument
essentially relies on the compact embeddings between some suitable
functional spaces of symmetric functions, which yields the existence of
nonnegative solutions of mountain-pass type, and the separation of the
corresponding mountain-pass levels from the energy levels associated to
radial solutions.\bigskip

 \smallskip


\end{abstract}

\section{Introduction and main result}

This paper is concerned with the following semilinear elliptic problem: 
\begin{equation}
\left\{ 
\begin{array}{l}
\begin{array}{ll}
\smallskip -\triangle u+\displaystyle\frac{A}{\left| x\right| ^{\alpha }}u=f\left(
u\right) ~ & \text{\textrm{in }}\mathbb{R}^{N},~N\geq 3 \\ 
\medskip u\geq 0 & \text{\textrm{in }}\mathbb{R}^{N}
\end{array}
\\ 
\begin{array}{l}
u\in H_{\alpha }^{1},\ u\neq 0
\end{array}
\end{array}
\right.  \tag*{$\left({\mathcal P}\right) $}
\end{equation}
where $A,\alpha >0$ are real constants, $f:\mathbb{R}\rightarrow \mathbb{R}$ is a
continuous nonlinearity satisfying $\left| f\left( s\right) \right| \leq
\left( \mathrm{const.}\right) s^{p-1}$ for some $p>2$ and all $s\geq 0$, and 
$H_{\alpha }^{1}:=D^{1,2}(\mathbb{R}^{N})\cap L^{2}(\mathbb{R}^{N},\left| x\right|
^{-\alpha }dx)$ is the natural energy space related to the equation. We will
deal with problem $\left( \mathcal{P}\right) $ in the weak sense, that is,
speaking about \emph{solutions} to $\left( \mathcal{P}\right) $ we will
always mean \emph{weak solutions}, i.e., functions $u\in H_{\alpha
}^{1}\setminus \left\{ 0\right\} $ such that $u\geq 0$ almost everywhere in $%
\mathbb{R}^{N}$ and 
\begin{equation}
\int_{\mathbb{R}^{N}}\nabla u\cdot \nabla v\,dx+\int_{\mathbb{R}^{N}}\displaystyle\frac{A}{%
\left| x\right| ^{\alpha }}uv\,dx=\int_{\mathbb{R}^{N}}f\left( u\right)
v\,dx\quad \text{for all }v\in H_{\alpha }^{1}.  \label{weak sol}
\end{equation}

As is well known, problems like $\left( \mathcal{P}\right) $ are models for
stationary states of reaction diffusion equations in population dynamics
(see e.g. \cite{Fife}) and arise in many branches of mathematical physics,
such as nonlinear optics, plasma physics, condensed matter physics and
cosmology (see e.g. \cite{BFmonograph,YangY}), where its nonnegative
solutions lead to special solutions (\textit{solitary waves} and \textit{%
solitons}) for several nonlinear field theories, like nonlinear Schr\"{o}%
dinger (or Gross-Pitaevskii) and Klein-Gordon equations. In this context, $%
\left( \mathcal{P}\right) $ is a prototype for problems exhibiting radial
potentials which are singular at the origin and/or vanishing at infinity
(sometimes called the \emph{zero mass case}; see e.g. \cite{Meder, BPR}).

\begin{figure}
\centering
\includegraphics{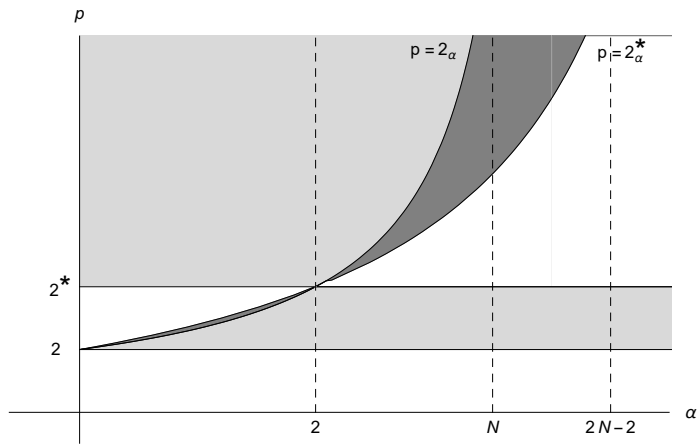} 
       \caption{Regions of nonexistence of solutions (light gray), and existence (white with $p>2$) and nonexistence (dark gray) of radial solutions.} 
								\label{pic}
\end{figure}

Though it can be considered of quite recent investigation, the study of
problem $\left( \mathcal{P}\right) $ has already some history, which
probably started in \cite{Terracini} and continued in \cite
{Co-Cr-Par,BRpow,Su-Wang-Will-2,Su-Wang-Will-p,BGRnonex, Catrina} (see \cite
{BGR} for a similar cylindrical problem). Currently, the problem of
existence and nonexistence of \emph{radial} solutions is essentially solved
in the pure-power case $f\left( u\right) =u^{p-1}$, where the results
obtained rest upon compatibility conditions between $\alpha $ and $p$. They
can be summarized as follows (for a chronological overview of these results
see \cite{BGRnonex}): the problem has a radial solution for $\left( \alpha
,p\right) =\left( 2,2^{*}\right) $ (\cite{Terracini}) and for all the pairs $%
\left( \alpha ,p\right) $ satisfying 
\begin{equation}
\left\{ 
\begin{array}{l}
0<\alpha <2 \\ 
2_{\alpha }^{*}<p<2^{*}
\end{array}
\right. \quad \text{or}\quad \left\{ 
\begin{array}{l}
2<\alpha <2N-2 \\ 
2^{*}<p<2_{\alpha }^{*}
\end{array}
\right. \quad \text{or}\quad \left\{ 
\begin{array}{l}
\alpha \geq 2N-2 \\ 
p>2^{*}
\end{array}
\right. ,\qquad 2_{\alpha }^{*}:=2\frac{2N-2+\alpha }{2N-2-\alpha }.
\label{sww}
\end{equation}
(\cite{Su-Wang-Will-p}), while it has no solution if 
\[
\left\{ 
\begin{array}{l}
0<\alpha <2 \\ 
p\notin \left( 2_{\alpha },2^{*}\right) 
\end{array}
\right. \quad \text{or}\quad \left\{ 
\begin{array}{l}
\alpha =2 \\ 
p\neq 2^{*}
\end{array}
\right. \quad \text{or}\quad \left\{ 
\begin{array}{l}
2<\alpha <N \\ 
p\notin \left( 2^{*},2_{\alpha }\right) 
\end{array}
\right. \quad \text{or}\quad \left\{ 
\begin{array}{l}
\alpha \geq N \\ 
p\leq 2^{*}
\end{array}
\right. ,\qquad 2_{\alpha }:=\frac{2N}{N-\alpha }
\]
(\cite{BRpow}) and no radial solution for both 
\[
\left\{ 
\begin{array}{l}
0<\alpha <2 \\ 
2_{\alpha }<p\leq 2_{\alpha }^{*}
\end{array}
\right. \quad \text{and}\quad \left\{ 
\begin{array}{l}
2<\alpha <2N-2 \\ 
2_{\alpha }^{*}\leq p<2_{\alpha }
\end{array}
\right. 
\]
(\cite{BGRnonex} and \cite{Catrina} respectively). As usual, $2^{*}:=2N/(N-2)
$ denotes the critical exponent for the Sobolev embedding in dimension $%
N\geq 3$. All these results are portrayed in the picture of the $\alpha p$%
-plane given in Fig.\ref{pic}, where nonexistence regions are shaded in
gray (nonexistence of radial solutions) and light gray (nonexistence of
solutions at all, which includes both the lines $p=2^{*}$ and $p=2_{\alpha }$
except for the pair $\left( \alpha ,p\right) =\left( 2,2^{*}\right) $),
whereas white color (of course above the line $p=2$) means existence of
radial solutions. 
As to nonradial solutions, the only result available is the one contained in 
\cite[Theorem 0.5]{Terracini}, where the author proves that problem $\left( 
\mathcal{P}\right) $ with $N\geq 4$, $\alpha =2$ and $f\left( u\right)
=u^{2^{*}-1}$ has at least a nonradial solution for every $A$ large enough.
This brings Catrina to say, in the introduction of his paper \cite{Catrina}:
``\textit{Two questions still remain: whether one can find non-radial
solutions in the case when radial solutions do not exist, or in the case
when radial solutions exist}''.

As concerns problem $\left( \mathcal{P}\right) $ with general nonlinearities
satisfying a power growth condition $\left| f\left( u\right) \right| \leq
\left( \mathrm{const.}\right) u^{p-1}$ for some $p>2$, the results of \cite
{Su-Wang-Will-p} also cover this case and ensures that, under some rather
standard additional assumptions on $f$ (precisely $(\mathbf{f}_{1})$ and $(%
\mathbf{f}_{2})$ below), $\left( \mathcal{P}\right) $ has a \emph{radial}
solution for all the pairs $\left( \alpha ,p\right) $ satisfying (\ref{sww})
again. To be precise, the authors only concern themselves with radial weak
solutions in the sense of the dual space of the radial subspace of $%
H_{\alpha }^{1}$ (where the energy functional of the problem is well defined
by the embeddings they prove), but the symmetric criticality type results of 
\cite{BGR-p2} actually apply, yielding solutions in the sense of our
definition (\ref{weak sol}). No results are known in the literature about
nonradial solutions.

This general lack of symmetry breaking results is the motivation of this
paper, where we prove that problem $\left( \mathcal{P}\right) $ has multiple
nonradial solutions as $A\rightarrow +\infty $ provided that $N\geq 4$, $%
\alpha \in (\frac{2}{N-1},2N-2)\setminus \left\{ 2\right\} $ and $f$ belongs
to a suitable class of nonlinearities satisfying a power growth condition.
We observe straight away that such a class of nonlinearities does not
unfortunately contain pure powers (which does not satisfy our assumption 
$(\mathbf{f}_{p_{1},p_{2}})$, where $p_{1}\neq p_{2}$).

The main assumptions characterizing our class of nonlinearities are the
following, where we denote $F\left( s\right) :=\int_{0}^{s}f\left( t\right)
dt$:

\begin{itemize}
\item[$(\mathbf{f}_{p_{1},p_{2}})$]  $\displaystyle \sup_{s>0}\frac{\left|
f\left( s\right) \right| }{\min \{s^{p_{1}-1},s^{p_{2}-1}\}}<+\infty $ for
some $2<p_{1}<2^{*}<p_{2}$

\item[$(\mathbf{f}_{1})$]  $\exists \theta >2$ such that $\theta F\left(
s\right) \leq f\left( s\right) s$ for all $s>0$

\item[$(\mathbf{f}_{2})$]  \text{$F\left( s\right) >0$ for }all $s>0$

\item[$(\mathbf{f}_{3})$]  the function $\text{$f$}$\text{$\left( s\right) /s
$} is strictly increasing on $\left( 0,+\infty \right) $

\item[$(\mathbf{f}_{4})$]  $\exists \mu >2$ such that the function $F$\text{$%
\left( s\right) /s^{\mu }$} is decreasing on $\left( 0,+\infty \right) .$
\end{itemize}

\noindent For $N\geq 3$, $\alpha \in \left(0,2N-2\right) $, $\alpha
\neq 2$, and $2<p_{1}<2^{*}<p_{2}$, we define 
\begin{equation}
\nu :=\nu _{N,\alpha ,p_{1},p_{2}}:=\left\{ 
\begin{array}{lll}
\left\lceil 2\min \left\{ \frac{N-1}{\alpha },\frac{N-2}{2-\alpha }\frac{%
2^{*}-p_{1}}{p_{1}-2}\right\} -2N\left( \frac{1}{\alpha }-\frac{1}{2}\right)
\right\rceil -1\medskip & \text{if } & 0<\alpha <2 \\ 
\left\lceil 2\min \left\{ \frac{N-1}{\alpha },\frac{N-2}{\alpha -2}\frac{%
p_{2}-2^{*}}{p_{2}-2}\right\} \right\rceil -1 & \text{if } & 2<\alpha <2N-2
\end{array}
\right.  \label{n_alfa}
\end{equation}
where $\left\lceil \cdot \right\rceil $ denotes the ceiling function (i.e., $%
\left\lceil x\right\rceil :=\min \left\{ n\in \mathbb{Z}:n\geq x\right\} $).

Our main result is the following theorem.

\begin{theorem}
\label{THM: main}Let $N\geq 4$ and $\alpha \in \left( 2/(N-1),2N-2\right) $, 
$\alpha \neq 2$. Let $f:\mathbb{R}\rightarrow \mathbb{R}$ be a continuous function
satisfying $(\mathbf{f}_{1})$-$(\mathbf{f}_{4})$. Assume that $(\mathbf{f}%
_{p_{1},p_{2}})$ holds with 
\begin{equation}
p_{1}<p_{1}^{*}:=2\frac{\alpha ^{2}(N-1)-2\alpha (N-1)+4N}{\alpha
^{2}(N-1)-2\alpha (N+1)+4N}\quad \text{or}\quad p_{2}>p_{2}^{*}:=2\frac{%
2N+2-\alpha }{2N-2-\alpha }  \label{p1p2}
\end{equation}
according as $\alpha \in \left( 2/(N-1),2\right) $ or $\alpha \in \left(
2,2N-2\right) $. Then there exists $A_{*}>0$ such that for every $A>A_{*}$
problem $\left( \mathcal{P}\right) $ has both a radial solution and $\nu $
different nonradial solutions. 
\end{theorem}

\begin{figure}
\centering
    \begin{subfigure}{0.49\textwidth}
\includegraphics[width=\textwidth,height=2.0in]
        {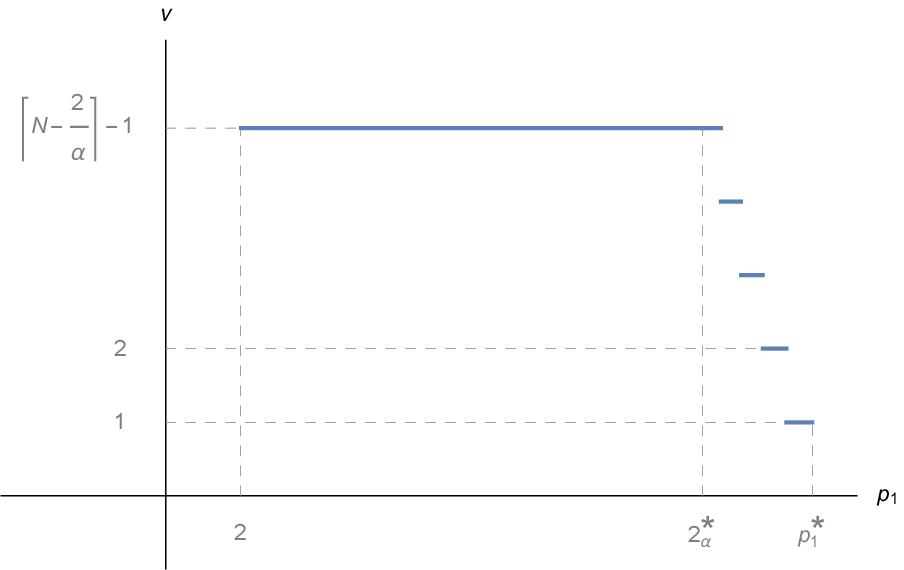}
       \caption{$\nu$ as a function of $p_1\in(2,p_1^*)$ for $N\geq 3$ and $\alpha\in(0,2)$ fixed.} 
								\label{p1}
\end{subfigure}
\begin{subfigure}{0.49\textwidth}
\includegraphics[width=\textwidth,height=2.0in]
        {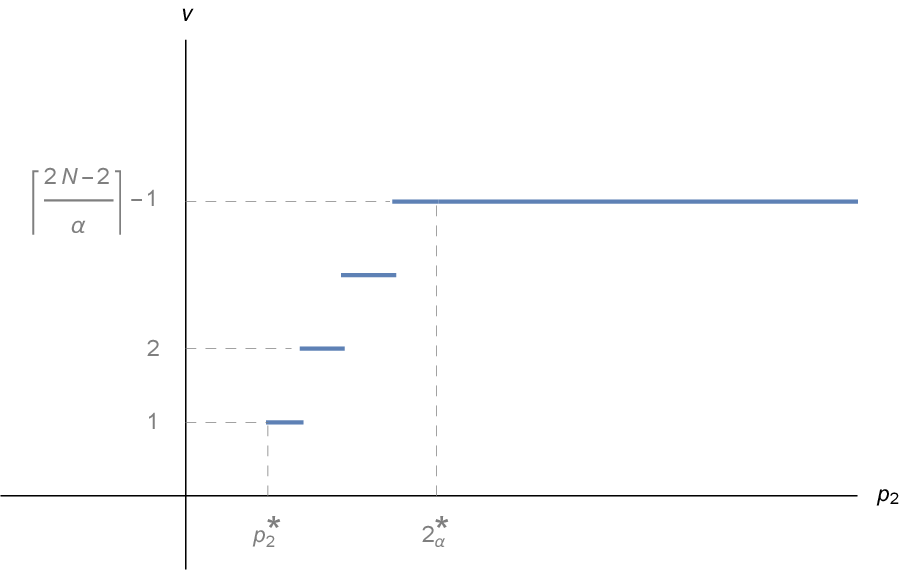}
        \caption{$\nu$ as a function of $p_2>p_2^*$ for $N\geq 3$ and $\alpha\in(2,2N-2)$ fixed.} 
        \label{p2}
\end{subfigure}
   \end{figure}

Some comments on Theorem \ref{THM: main} are in order. First of all, under
the assumptions of the theorem, $\nu $ is positive (see Lemma \ref{LEM:
conti} below), so that at least one nonradial solution actually exists. On
the other hand, it is easy to check that, for every fixed $N$ and $\alpha $,
the behaviour of $\nu $ as a function of $p_{1}$ and $p_{2}$ is the one
portrayed in Figs. (\ref{p1}) and (\ref{p2}), whence one sees that the number 
$\nu $ of nonradial solutions may assume every natural value (as $N\rightarrow \infty $).

As to assumptions (\ref{p1p2}) and $(\mathbf{f}_{p_{1},p_{2}})$, it is worth
observing that $\alpha \in \left( 0,2\right) $ implies $2<p_{1}^{*}<2^{*}$,
while $2<\alpha <2N-2$ implies $p_{2}^{*}>2^{*}$, so that (\ref{p1p2}) and $(%
\mathbf{f}_{p_{1},p_{2}})$ are consistent with each other. Assumption $(%
\mathbf{f}_{p_{1},p_{2}})$ is the so-called \textit{double-power growth
condition} and seems to be typical in nonlinear problems with potentials
vanishing at infinity (see e.g. \cite
{Be-Gr-Mic,Be-Gr-Mic-2,Azz-Pomp,G-R,BGR-p2,BGR-p,BPR,Ghim-Mic,
Meder,Benci-Mic-exterior,BBR-2}). It obviously implies the single-power growth
condition $\left| f\left( s\right) \right| \leq \left( \mathrm{const.}%
\right) s^{p-1}$ for all $p\in [p_{1},p_{2}]$ and $s\geq 0$, but it is actually more
stringent than that, since it requires $p_{1}\neq p_{2}$. We finally observe
that $(\mathbf{f}_{p_{1},p_{2}})$ still remains true if one raises $p_{1}$
and lowers $p_{2}$, but this decreases $\nu $ (see Figs. (\ref{p1}) and (\ref{p2})) and
therefore it is convenient to apply Theorem \ref{THM: main} with $p_{1}$ as
small as possible $p_{2}$ as large as possible (which is also consistent
with assumption (\ref{p1p2})).

The plan of the paper is the following. In Section \ref{SEC:prelim} we define the
variational setting and introduce the argument we will use in the proof of
Theorem \ref{THM: main}, which will be completed in Section \ref{SEC: pf}.
Observe that we cannot use the technique used in \cite{Terracini}, where the
homogeneity of the nonlinearity is exploited and a nonradial solution is
obtained as a global minimizer of the Sobolev-type quotient associated to
the problem. Our argument, instead, essentially relies on the following two
main elements: (i) the compact embeddings between some suitable functional
spaces of symmetric functions, which yields the existence of $\nu $
different solutions of mountain-pass type; (ii) the separation of the
corresponding mountain-pass levels from the energy levels associated to
radial solutions. Sections \ref{SEC: mr} and \ref{SEC: cs} are devoted to
the estimation of these levels in order to separate them.
As a conclusion, we get $\nu$ nonradial solutions on which the energy functional 
of the equation has a lower value than the energy levels of radial solutions.

We end this introductory section by giving some examples of nonlinearities
to which Theorem \ref{THM: main} applies and collecting some notations we
use throughout the paper.

\begin{example}
Let $N\geq 4$ and $\alpha \in \left( 2/(N-1),2N-2\right) $, $\alpha \neq 2$,
and let $2<p_{1}<2^{*}<p_{2}$ be such that (\ref{p1p2}) holds. The most
obvious nonlinearity to which Theorem \ref{THM: main} applies is $f\left(
s\right) =\min \{\left| s\right| ^{p_{1}-1},\left| s\right| ^{p_{2}-1}\}$,
which satisfies $(\mathbf{f}_{1})$ and $(\mathbf{f}_{4})$ for $\theta =p_{1}$
and any $\mu >p_{2}$. Other simple examples are 
\[
f\left( s\right) =\frac{\left| s\right| ^{p_{2}-1}}{1+\left| s\right|
^{p_{2}-p_{1}}},\quad f\left( s\right) =\frac{d}{ds}\left( \frac{\left|
s\right| ^{p_{2}}}{1+\left| s\right| ^{p_{2}-p_{1}}}\right) ,
\]
both of which satisfy $(\mathbf{f}_{1})$ with $\theta =p_{1}$. In the latter
case, $(\mathbf{f}_{4})$ clearly holds for any $\mu >p_{2}$. We leave to the
reader to check that $(\mathbf{f}_{4})$ also holds in the former case, for $%
\mu $ large enough.
\end{example}

\noindent \textbf{Notations.}\smallskip

\noindent $\bullet $ $\sigma _{d}$ denotes the $\left( d-1\right) $%
-dimensional measure of the unit sphere of $\mathbb{R}^{d}$.\smallskip

\noindent $\bullet $ $C_{c}^{\infty }\left( \Omega \right) $ is the space of
the infinitely differentiable real functions with compact support in the
open set $\Omega \subseteq \mathbb{R}^{d}$.\smallskip

\noindent $\bullet $ $D^{1,2}(\mathbb{R}^{N})=\{u\in L^{2^{*}}(\mathbb{R}%
^{N}):\nabla u\in L^{2}(\mathbb{R}^{N})\}$ is the usual Sobolev space, which
identifies with the completion of $C_{c}^{\infty }(\mathbb{R}^{N})$ with
respect to the norm of the gradient.

\section{Preliminaries \label{SEC:prelim}}

Let $N\geq 3$ and $A,\alpha >0$. Let $f:\mathbb{R}\rightarrow \mathbb{R}$ be a
continuous function satisfying $(\mathbf{f}_{p_{1},p_{2}})$, $(\mathbf{f}%
_{1})$, $(\mathbf{f}_{2})$. In this section we define the functional setting
and introduce the argument we will use in proving Theorem \ref{THM: main}.

As already mentioned in the introduction, we define the Hilbert space 
\[
H_{\alpha }^{1}:=\left\{ u\in D^{1,2}(\mathbb{R}^{N}):\int_{\mathbb{R}^{N}}\frac{%
u^{2}}{\left| x\right| ^{\alpha }}dx<\infty \right\} , 
\]
which we endow with the following scalar product and related norm: 
\[
\left( u,v\right) _{A}:=\int_{\mathbb{R}^{N}}\left( \nabla u\cdot \nabla v+%
\frac{A}{\left| x\right| ^{\alpha }}uv\right) dx,\quad \left\| u\right\|
_{A}^{2}:=\int_{\mathbb{R}^{N}}\left( \left| \nabla u\right| ^{2}+\frac{A}{%
\left| x\right| ^{\alpha }}u^{2}\right) dx. 
\]
Here and in the rest of the paper, $D^{1,2}(\mathbb{R}^{N})=\{u\in L^{2^{*}}(%
\mathbb{R}^{N}):\nabla u\in L^{2}(\mathbb{R}^{N})\}$ is the usual Sobolev
space, which identifies with the completion of $C_{\mathrm{c}}^{\infty }(%
\mathbb{R}^{N})$ with respect to the norm of the gradient. Of course, the
embedding $H_{\alpha }^{1}\hookrightarrow D^{1,2}(\mathbb{R}^{N})$ is
continuous. 

Given any integer $K$ such that $1\leq K\leq N-1$, we write every $x\in \mathbb{%
R}^{N}$ as $x=\left( y,z\right) \in \mathbb{R}^{K}\times \mathbb{R}^{N-K}$ and in
the space $H_{\alpha }^{1}$ we consider the following closed subspaces of
symmetric functions: 
\[
H_{\mathrm{r}}:=\left\{ u\in H:u\left( x\right) =u\left( \left| x\right|
\right) \right\} \quad \text{and}\quad H_{K}:=\left\{ u\in H:u\left(
x\right) =u\left( y,z\right) =u\left( \left| y\right| ,\left| z\right|
\right) \right\} . 
\]
Of course $u\left( y,z\right) =u\left( \left| y\right| ,\left| z\right|
\right) $ naturally means that $u\left( y,z\right) =u\left(
S_{1}y,S_{2}z\right) $ for all isometries $S_{1}$ and $S_{2}$ of $\mathbb{R}%
^{K} $ and $\mathbb{R}^{N-K}$ respectively. Similarly for $u\left( x\right)
=u\left( \left| x\right| \right) $. Note that $H_{\mathrm{r}}\subset H_{K}$
for every $K$, since $\left| x\right| ^{2}=\left| y\right| ^{2}+\left|
z\right| ^{2}$. The next lemma better clarifies the relation between the
spaces $H_{K}$ and $H_{\mathrm{r}}$.

\begin{lemma}
\label{LEM: Hk}Let $1\leq K_{1}<K_{2}\leq N-1$. Then $H_{K_{1}}\cap
H_{K_{2}}=H_{\mathrm{r}}$.
\end{lemma}

\proof
The proof is essentially an adaptation of the one of \cite[Lemma 3.3]{Li-90}%
. We will denote by $\left( y_{1},z_{1}\right) $ the decomposition in $\mathbb{R%
}^{K_{1}}\times \mathbb{R}^{N-K_{1}}$ of any $x\in \mathbb{R}^{N}$, and by $\left(
y_{2},z_{2}\right) $ its decomposition in $\mathbb{R}^{K_{2}}\times \mathbb{R}%
^{N-K_{2}}$. 
Let $u\in H_{K_{1}}\cap H_{K_{2}}$ and for every $s,t\geq 0$ define 
\[
\tilde{u}\left( s,t\right) :=u\left( s,0,...,0,t\right) . 
\]
Then we clearly have $u\left( x\right) =\tilde{u}\left( \left| y_{1}\right|
,\left| z_{1}\right| \right) =\tilde{u}\left( \left| y_{2}\right| ,\left|
z_{2}\right| \right) $ for every $x=\left( y_{1},z_{1}\right) =\left(
y_{2},z_{2}\right) \in \mathbb{R}^{N}$. Let $x=\left( y_{1},z_{1}\right) \in 
\mathbb{R}^{K_{1}}\times \mathbb{R}^{N-K_{1}}$ and $x^{\prime }=\left(
y_{2}^{\prime },z_{2}^{\prime }\right) \in \mathbb{R}^{K_{2}}\times \mathbb{R}%
^{N-K_{2}}$ be such that $\left| x\right| =\left| x^{\prime }\right| $,
i.e., $\left| y_{1}\right| ^{2}+\left| z_{1}\right| ^{2}=\left|
y_{2}^{\prime }\right| ^{2}+\left| z_{2}^{\prime }\right| ^{2}$. Suppose
that $\left| y_{1}\right| \leq \left| y_{2}^{\prime }\right| $ and define $%
x^{\prime \prime }\in \mathbb{R}^{N}$ by setting 
\[
x^{\prime \prime }:=\left( \left| y_{1}\right| ,0,...,0,\sqrt{\left|
y_{2}^{\prime }\right| ^{2}-\left| y_{1}\right| ^{2}},0...,0,\left|
z_{2}^{\prime }\right| ,0,...,0\right) 
\]
where the first block of zeros has $K_{1}-1$ zeros, the second $%
K_{2}-K_{1}-1 $, and the third $N-K_{2}-1$. Then 
\[
u\left( x^{\prime \prime }\right) =\tilde{u}\left( \left| y_{1}^{\prime
\prime }\right| ,\left| z_{1}^{\prime \prime }\right| \right) =\tilde{u}%
\left( \left| y_{1}\right| ,\sqrt{\left| y_{2}^{\prime }\right| ^{2}-\left|
y_{1}\right| ^{2}+\left| z_{2}^{\prime }\right| ^{2}}\right) =\tilde{u}%
\left( \left| y_{1}\right| ,\left| z_{1}\right| \right) 
\]
and 
\[
u\left( x^{\prime \prime }\right) =\tilde{u}\left( \left| y_{2}^{\prime
\prime }\right| ,\left| z_{2}^{\prime \prime }\right| \right) =\tilde{u}%
\left( \sqrt{\left| y_{1}\right| ^{2}+\left| y_{2}^{\prime }\right|
^{2}-\left| y_{1}\right| ^{2}},\left| z_{2}^{\prime }\right| \right) =\tilde{%
u}\left( \left| y_{2}^{\prime }\right| ,\left| z_{2}^{\prime }\right|
\right) , 
\]
which implies $\tilde{u}\left( \left| y_{1}\right| ,\left| z_{1}\right|
\right) =\tilde{u}\left( \left| y_{2}^{\prime }\right| ,\left| z_{2}^{\prime
}\right| \right) $, i.e., $u\left( x\right) =u\left( x^{\prime }\right) $.
If $\left| y_{1}\right| >\left| y_{2}^{\prime }\right| $ we repeat the
argument with $x^{\prime \prime }$ defined by 
\[
x^{\prime \prime }:=\left( \left| y_{2}^{\prime }\right| ,0,...,0,\sqrt{%
\left| y_{1}\right| ^{2}-\left| y_{2}^{\prime }\right| ^{2}},0...,0,\left|
z_{1}\right| ,0,...,0\right) 
\]
(with the same lenghts of the blocks) and get the same result. Hence $\left|
x\right| =\left| x^{\prime }\right| $ implies $u\left( x\right) =u\left(
x^{\prime }\right) $, i.e., $u\in H_{\mathrm{r}}$.%
\endproof

We modify the function $f$ by setting $f\left( s\right) =0$ for all $s<0$
and, with a slight abuse of notation, we still denote by $f$ the modified
function. Then by $(\mathbf{f}_{p_{1},p_{2}})$ there exist $M_{1},M_{2}>0$
such that 
\[
\left| f\left( s\right) \right| \leq M_{1}\min \{\left| s\right|
^{p_{1}-1},\left| s\right| ^{p_{2}-1}\}\quad \text{and}\quad \left| F\left(
s\right) \right| \leq M_{2}\min \{\left| s\right| ^{p_{1}},\left| s\right|
^{p_{2}}\}\quad \text{for all }s\in \mathbb{R}, 
\]
which yields in particular 
\begin{equation}
\left| f\left( s\right) \right| \leq M_{1}\left| s\right| ^{p-1}\quad \text{%
and}\quad \left| F\left( s\right) \right| \leq M_{2}\left| s\right|
^{p}\quad \text{for all }p\in \left[ p_{1},p_{2}\right] \text{ and }s\in 
\mathbb{R}.  \label{growth p}
\end{equation}
By the continuous embeddings $H_{\alpha }^{1}\hookrightarrow D^{1,2}(\mathbb{R}%
^{N})\hookrightarrow L^{2^{*}}(\mathbb{R}^{N})$, one checks (see for example 
\cite{Kuz-Poho}) that condition (\ref{growth p}) with $p=2^{*}$ implies that
the energy functional associated to the equation of $\left( \mathcal{P}\right) $, i.e.,
\begin{equation}
I\left( u\right) :=\frac{1}{2}\left\| u\right\| _{A}^{2}-\int_{\mathbb{R}%
^{N}}F\left( u\right) dx,
  \label{I(u)=}
\end{equation}
is of class $C^{1}$ on $H_{\alpha }^{1}$ and has Fr\'{e}chet derivative $%
I^{\prime }\left( u\right) $ at any $u\in H_{\alpha }^{1}$ given by 
\begin{equation}
I^{\prime }\left( u\right) v=\left( u,v\right) _{A}-\int_{\mathbb{R}%
^{N}}f\left( u\right) v\,dx\quad \text{for all }v\in H_{\alpha }^{1}.
\label{I'(u)h=}
\end{equation}
This yields that critical points of $I:H_{\alpha }^{1}\rightarrow \mathbb{R}$
satisfy (\ref{weak sol}). A standard argument shows that such critical
points are nonnegative (see the proof of Theorem \ref{THM: main} in Section 
\ref{SEC: pf}) and therefore we conclude that nonzero critical points of $I$
are weak solutions to problem $\left( \mathcal{P}\right) $. %

Accordingly, our argument in proving Theorem \ref{THM: main} will be
essentially the following. The existence of a critical point for the
restriction $I_{\mid H_{\mathrm{r}}}$ readily follows from the results of 
\cite{Su-Wang-Will-p}. By exploiting the compact embeddings of \cite
{Azz-Pomp} and the results of \cite{BPR} about Nemytski\u{\i} operators on
the sum of Lebesgue spaces, we will show in Section \ref{SEC: pf} that $%
I_{\mid H_{K}}$ has a nonzero critical point $u_{K}$ for every $2\leq K\leq
N-2$. Thanks to the classical Palais' Principle of Symmetric Criticality 
\cite{Palais}, all these critical points are also critical points of $I$,
and thus weak solutions to $\left( \mathcal{P}\right) $. Hence Theorem \ref
{THM: main} is proved is we show that $u_{K}\notin H_{\mathrm{r}}$ for every 
$K$, which also implies $u_{K_{1}}\neq u_{K_{2}}$ for $K_{1}\neq K_{2}$ by
Lemma \ref{LEM: Hk}. This will achieved by showing that the critical levels $%
I\left( u_{K}\right) $ are lower than all the nonzero critical levels of $%
I_{\mid H_{\mathrm{r}}}$. The starting points in proving this are the
following lemmas.

\begin{lemma}
\label{LEM: max}For every $u\in H_{\mathrm{r}}\setminus \left\{ 0\right\} $, 
$u\geq 0$, there exists $t_{u}>0$ such that $I\left( t_{u}u\right)
=\max_{t\geq 0}I\left( tu\right) $.
\end{lemma}

\proof
Since $u\geq 0$ and $u\neq 0$, we can fix $\delta >0$ such that the set $%
\left\{ x\in \mathbb{R}^{N}:u\geq \delta \right\} $ has positive measure. From
assumptions $(\mathbf{f}_{1})$ and $(\mathbf{f}_{2})$, we deduce that there
exists a constant $C>0$ such that $F\left( s\right) \geq Cs^{\theta }$ for
all $s\geq \theta $. Then for every $t>1$ one has 
\begin{eqnarray*}
\int_{\mathbb{R}^{N}}F\left( tu\right) dx &=&\int_{\left\{ x\in \mathbb{R}%
^{N}:tu\geq \delta \right\} }F\left( u\right) dx+\int_{\left\{ x\in \mathbb{R}%
^{N}:tu<\delta \right\} }F\left( u\right) dx\geq \int_{\left\{ x\in \mathbb{R}%
^{N}:tu\geq \delta \right\} }F\left( u\right) dx \\
&\geq &Ct^{\theta }\int_{\left\{ x\in \mathbb{R}^{N}:tu\geq \delta \right\}
}u^{\theta }dx\geq Ct^{\theta }\int_{\left\{ x\in \mathbb{R}^{N}:u\geq \delta
\right\} }u^{\theta }dx
\end{eqnarray*}
and therefore 
\[
I\left( tu\right) \leq \frac{1}{2}t^{2}\left\| u\right\| _{A}^{2}-Ct^{\theta
}\int_{\left\{ x\in \mathbb{R}^{N}:u\geq \delta \right\} }u^{\theta
}dx\rightarrow -\infty \quad \text{as }t\rightarrow +\infty 
\]
since $\theta >2$. As $I\left( 0\right) =0$, this gives the result.%
\endproof

According to Lemma \ref{LEM: max}, define 
\[
m_{A}:=\inf_{u\in H_{\mathrm{r}}\setminus \left\{ 0\right\} ,\,u\geq
0}\max_{t\geq 0}I\left( tu\right) . 
\]

\begin{lemma}
\label{LEM: mr}Assume $(\mathbf{f}_{3})$ and let $u\in H_{\mathrm{r}%
}\setminus \left\{ 0\right\} $. If $u$ is a critical point for $I$, then 
\[
I\left( u\right) \geq m_{A}.
\]
\end{lemma}

\proof
As already observed, if $u$ is a critical point for $I$ then $u$ is
nonnegative. We now prove that $I\left( u\right) =\max_{t\geq 0}I\left(
tu\right) $, which obviously yields the result. For $t\geq 0$ define 
\[
g\left( t\right) :=I\left( tu\right) =\frac{1}{2}t^{2}\left\| u\right\|
_{A}^{2}-\int_{\mathbb{R}^{N}}F\left( tu\right) dx. 
\]
As $u$ is a critical point for $I$, we readily have that $t=1$ is a critical
point for $g$. Indeed, $g^{\prime }\left( t\right) =I^{\prime }\left(
tu\right) u$ and thus $g^{\prime }\left( 1\right) =I^{\prime }\left(
u\right) u=0$. We now show that, on the other hand, $g$ has at most one
critical point in $\left( 0,+\infty \right) $. We have $g^{\prime }\left(
t\right) =0$ if and only if $I^{\prime }\left( tu\right) u=0$, i.e. 
\[
I^{\prime }\left( tu\right) u=t\left\| u\right\| _{A}^{2}-\int_{\mathbb{R}%
^{N}}f\left( tu\right) u\,dx=0. 
\]
So, if $0<t_{1}<t_{2}$ are critical points for $g$, one has 
\[
\left\| u\right\| _{A}^{2}=\frac{1}{t_{1}}\int_{\mathbb{R}^{N}}f\left(
t_{1}u\right) u\,dx=\frac{1}{t_{2}}\int_{\mathbb{R}^{N}}f\left( t_{2}u\right)
u\,dx, 
\]
which implies 
\begin{equation}
\int_{E_{u}}\left( \frac{f\left( t_{2}u\right) }{t_{2}u}-\frac{f\left(
t_{1}u\right) }{t_{1}u}\right) u^{2}dx=0  \label{f3}
\end{equation}
where $E_{u}:=\left\{ x\in \mathbb{R}^{N}:u>0\right\} $. Since the integrand in
(\ref{f3}) is nonnegative by assumption $(\mathbf{f}_{3})$, we have that $%
f\left( t_{2}u\right) /\left( t_{2}u\right) -f\left( t_{1}u\right) /\left(
t_{1}u\right) =0$ almost everywhere on $E_{u}$. Since $E_{u}$ has positive
measure (because $0\neq u\geq 0$ ), this implies $t_{1}=t_{2}$, again by
assumption $(\mathbf{f}_{3})$. As a conclusion, according to Lemma \ref{LEM:
max}, we deduce that $t_{u}=1$ and the claim ensues.%
\endproof

\begin{lemma}
\label{LEM: mp}There exist $R>0$ such that for every $1\leq K\leq N-1$ one
has 
\begin{equation}   \label{mp0}
\inf_{u\in H_{K},\left\| u\right\| _{A}\leq R}I\left( u\right) =0\quad \text{%
and\quad }\inf_{u\in H_{K},\left\| u\right\| _{A}=R}I\left( u\right) >0.
\end{equation}
\end{lemma}

\proof
The claim readily follows from (\ref{growth p}) with $p=2^{*}$ and the
continuous embeddings $H_{K}\hookrightarrow H_{\alpha }^{1}\hookrightarrow
D^{1,2}(\mathbb{R}^{N})\hookrightarrow L^{2^{*}}(\mathbb{R}^{N})$, which imply
that there exists a constant $C>0$ such that 
$I\left( u\right) \geq \left\| u\right\| _{A}^{2}/2-C\left\|u\right\| _{A}^{2^{*}}$
for all $u\in H_{K}$.%
\endproof

In Section \ref{SEC: cs} we will see that $I_{\mid H_{K}}$ takes negative
values by choosing a suitable $\overline{u}_{K}\in H_{K}$ such that $I\left( 
\overline{u}_{K}\right) <0$. This implies $\left\| \overline{u}_{K}\right\|
_{A}>R$ by \eqref{mp0} 
and therefore the functional $I_{\mid H_{K}}$
has a mountain-pass geometry. In Section \ref{SEC: pf} we will see that it
also satisfies the Palais-Smale condition for $2\leq K\leq N-2$, so that it
admits a (nonnegative) critical point $u_{K}$ at the mountain-pass level 
\begin{equation}
c_{A,K}:=\inf_{\gamma \in \Gamma }\max_{t\in \left[ 0,1\right] }I\left(
\gamma \left( t\right) \right) >0\quad \text{where\quad }\Gamma :=\left\{
\gamma \in C\left( \left[ 0,1\right] ;H_{K}\right) :\gamma \left( 0\right)
=0,\,\gamma \left( 1\right) =\overline{u}_{K}\right\} .  \label{cs:=}
\end{equation}
With a view to obtaining the separation inequality $c_{A,K}<m_{A}$, Sections 
\ref{SEC: mr} and \ref{SEC: cs} are devoted to estimating $m_{A}$ and $%
c_{A,K}$. 

{\allowdisplaybreaks
\section{Estimate of $m_{A}$ \label{SEC: mr}}

Let $N\geq 3$ and $\alpha ,A>0$. Let $f:\mathbb{R}\rightarrow \mathbb{R}$ be a
continuous function satisfying $(\mathbf{f}_{p_{1},p_{2}})$, $(\mathbf{f}%
_{1})$, $(\mathbf{f}_{2})$. This section is devoted to deriving the estimate
of $m_{A}$ given in Proposition \ref{PROP: mr>} below, which relies on the
following radial lemma (see also \cite[Appendix]{BGR-p2} and 
\cite[Lemmas 4 and 5]{Su-Wang-Will-p} for similar results).

\begin{lemma}
\label{LEM: radial}Every $u\in H_{\mathrm{r}}$ satisfies 
\[
\left| u\left( x\right) \right| \leq \frac{\sqrt{2/\sigma _{N}}}{A^{1/4}}%
\frac{\left\| u\right\| _{A}}{\left| x\right| ^{\frac{2N-2-\alpha }{4}}}%
\quad \text{almost everywhere in }\mathbb{R}^{N}
\]
(recall that $\sigma _{N}$ denotes the $\left( N-1\right) $-dimensional
measure of the unit sphere of $\mathbb{R}^{N}$).
\end{lemma}

\proof
Let $u\in H_{\mathrm{r}}$ and let $\tilde{u}:\left( 0,+\infty \right)
\rightarrow \mathbb{R}$ be continuous and such that $u\left( x\right) =\tilde{u}%
\left( \left| x\right| \right) $ for almost every $x\in \mathbb{R}^{N}$. Set 
\[
v\left( r\right) :=r^{N-1-\alpha /2}\tilde{u}\left( r\right) ^{2}\quad \text{%
for all }r>0. 
\]
By \cite[Lemma 27]{BGR-p2} we have that $\tilde{u}\in W^{1,1}\left( \left(
a,b\right) \right) $ for every $0<a<b<+\infty $, whence $v\in W^{1,1}\left(
\left( a,b\right) \right) $ and 
\[
v\left( b\right) -v\left( a\right) =\int_{a}^{b}v^{\prime }\left( r\right)
dr\,. 
\]
Moreover, for almost every $r\in \left( a,b\right) $ one has 
\begin{equation}
v^{\prime }\left( r\right) =\left( N-1-\frac{\alpha }{2}\right)
r^{N-2-\alpha /2}\tilde{u}\left( r\right) ^{2}+2r^{N-1-\alpha /2}\tilde{u}%
\left( r\right) \tilde{u}^{\prime }\left( r\right) .  \label{v'=}
\end{equation}
If $\alpha <2N-2$, this implies $v^{\prime }\left( r\right) \geq
2r^{N-1-\alpha /2}\tilde{u}\left( r\right) \tilde{u}^{\prime }\left(
r\right) $ and therefore 
\begin{eqnarray}
v\left( a\right) &\leq &v\left( b\right) -\int_{a}^{b}v^{\prime }\left(
r\right) dr\leq v\left( b\right) -\int_{a}^{b}2r^{N-1-\frac{\alpha }{2}}%
\tilde{u}\left( r\right) \tilde{u}^{\prime }\left( r\right) dr  \label{v(a)<}
\\
&\leq &v\left( b\right) +2\int_{a}^{b}r^{N-1-\frac{\alpha }{2}}\left| \tilde{%
u}\left( r\right) \right| \left| \tilde{u}^{\prime }\left( r\right) \right|
dr=v\left( b\right) +2\int_{a}^{b}r^{\frac{N-1}{2}}\left| \tilde{u}^{\prime
}\left( r\right) \right| r^{\frac{N-1}{2}-\frac{\alpha }{2}}\left| \tilde{u}%
\left( r\right) \right| dr  \nonumber \\
&\leq &v\left( b\right) +2\left( \int_{0}^{+\infty }r^{N-1}\left| \tilde{u}%
^{\prime }\left( r\right) \right| ^{2}dr\right) ^{1/2}\left(
\int_{0}^{+\infty }r^{N-1-\alpha }\left| \tilde{u}\left( r\right) \right|
^{2}dr\right) ^{1/2}  \nonumber \\
&\leq &v\left( b\right) +\frac{2}{\sigma _{N}\sqrt{A}}\left( \int_{\mathbb{R}%
^{N}}\left| \nabla u\right| ^{2}dx\right) ^{1/2}\left( \int_{\mathbb{R}^{N}}%
\frac{Au^{2}}{\left| x\right| ^{\alpha }}dx\right) ^{1/2}\leq \frac{2}{%
\sigma _{N}\sqrt{A}}\left\| u\right\| _{A}^{2}  \nonumber
\end{eqnarray}
If $\alpha \geq 2N-2$, then (\ref{v'=}) gives $v^{\prime }\left( r\right)
\leq 2r^{N-1-\alpha /2}\tilde{u}\left( r\right) \tilde{u}^{\prime }\left(
r\right) $ and thus we have 
\begin{equation}
v\left( b\right) \leq v\left( a\right) +\int_{a}^{b}2r^{N-1-\alpha /2}\tilde{%
u}\left( r\right) \tilde{u}^{\prime }\left( r\right) dr\leq v\left( a\right)
+\frac{2}{\sigma _{N}\sqrt{A}}\left\| u\right\| _{A}^{2}  \label{v(b)<}
\end{equation}
as before. Now observe that there exist $0<a_{n}\rightarrow 0$ and $%
b_{n}\rightarrow +\infty $ such that $v\left( a_{n}\right) \rightarrow 0$
and $v\left( b_{n}\right) \rightarrow 0$. Indeed, if $l:=\liminf_{r%
\rightarrow 0^{+}}v\left( r\right) >0$, then for every $r$ smaller than some
suitable $r_{0}>0$ one has $\left| \tilde{u}\left( r\right) \right| \geq 
\sqrt{l/2}\,r^{-(N-1-\alpha /2)/2}$ and therefore one of the following
contradictions ensues: 
\[
\int_{\mathbb{R}^{N}}\frac{u^{2}}{\left| x\right| ^{\alpha }}dx\geq \frac{l}{2}%
\int_{B_{r_{0}}}\frac{1}{\left| x\right| ^{N-1+\alpha /2}}dx=+\infty \quad 
\text{if }\alpha \geq 2 
\]
or 
\[
\int_{\mathbb{R}^{N}}\left| u\right| ^{2^{*}}dx\geq \left( \frac{l}{2}\right)
^{2^{*}/2}\int_{B_{r_{0}}}\frac{1}{\left| x\right| ^{\frac{N-1-\alpha /2}{N-2%
}N}}dx=+\infty \quad \text{if }\alpha \leq 2. 
\]
Similary, if $\liminf_{r\rightarrow +\infty }v\left( r\right) >0$, one
obtains $\int_{\mathbb{R}^{N}}\frac{u^{2}}{\left| x\right| ^{\alpha }}%
dx=+\infty $ if $\alpha \leq 2$ and $\int_{\mathbb{R}^{N}}\left| u\right|
^{2^{*}}dx=+\infty $ if $\alpha \geq 2$. Hence the claim follows by letting $%
n\rightarrow \infty $ in (\ref{v(a)<}) with $a=r$ and $b=b_{n}$, and in (\ref
{v(b)<}) with $a=a_{n}$ and $b=r$.%
\endproof

We can now prove our estimate for $m_{A}$.

\begin{proposition}
\label{PROP: mr>}Assume $0<\alpha <2N-2$, $\alpha \neq 2$, and let $p=\max
\{2_{\alpha }^{*},p_{1}\}$ or $p=\min \{2_{\alpha }^{*},p_{2}\}$ according
as $0<\alpha <2$ or $2<\alpha <2N-2$. Then there exists a constant $C_{0}>0$
independent from $A$ such that 
\[
m_{A}\geq C_{0}A^{\frac{N-2}{\alpha -2}\frac{p-2^{*}}{p-2}}.
\]
\end{proposition}

\proof
Let $u\in H_{\mathrm{r}}\setminus \left\{ 0\right\} $. By Lemma \ref{LEM:
radial}, we have 
\begin{eqnarray*}
\int_{\mathbb{R}^{N}}\left| u\right| ^{2_{\alpha }^{*}}dx &=&\int_{\mathbb{R}%
^{N}}\left| u\right| ^{2_{\alpha }^{*}-2}u^{2}dx\leq \frac{\left( 2/\sigma
_{N}\right) ^{\left( 2_{\alpha }^{*}-2\right) /2}}{A^{\left( 2_{\alpha
}^{*}-2\right) /4}}\left\| u\right\| _{A}^{2_{\alpha }^{*}-2}\int_{\mathbb{R}%
^{N}}\frac{u^{2}}{\left| x\right| ^{\frac{2N-2-\alpha }{4}\left( 2_{\alpha
}^{*}-2\right) }}dx \\
&=&\frac{\left( 2/\sigma _{N}\right) ^{\frac{2\alpha }{2N-2-\alpha }}}{A^{%
\frac{\alpha }{2N-2-\alpha }}}\frac{\left\| u\right\| _{A}^{2_{\alpha
}^{*}-2}}{A}\int_{\mathbb{R}^{N}}\frac{Au^{2}}{\left| x\right| ^{\alpha }}%
dx\leq \frac{\left( 2/\sigma _{N}\right) ^{\frac{2\alpha }{2N-2-\alpha }}}{%
A^{\frac{2N-2}{2N-2-\alpha }}}\left\| u\right\| _{A}^{2_{\alpha }^{*}}
\end{eqnarray*}
since $2_{\alpha }^{*}-2=4\alpha /\left( 2N-2-\alpha \right) $. On the other
hand, one has 
\[
\int_{\mathbb{R}^{N}}\left| u\right| ^{2^{*}}dx\leq S_{N}^{2^{*}}\left( \int_{%
\mathbb{R}^{N}}\left| \nabla u\right| ^{2}dx\right) ^{2^{*}/2}\leq
S_{N}^{2^{*}}\left\| u\right\| _{A}^{2^{*}} 
\]
where $S_{N}$ denotes the Sobolev constant in dimension $N$. Then, both for $%
p=\max \{2_{\alpha }^{*},p_{1}\}<2^{*}$ and $p=\min \{2_{\alpha
}^{*},p_{2}\}>2^{*}$, we can argue by interpolation: there exists $\lambda
\in \left( 0,1\right) $ such that $p=\lambda 2^{*}+\left( 1-\lambda \right)
2_{\alpha }^{*}$ and by H\"{o}lder inequality we get 
\[
\int_{\mathbb{R}^{N}}\left| u\right| ^{p}dx=\int_{\mathbb{R}^{N}}\left| u\right|
^{\lambda 2^{*}+\left( 1-\lambda \right) 2_{\alpha }^{*}}dx\leq \left( \int_{%
\mathbb{R}^{N}}\left| u\right| ^{2^{*}}dx\right) ^{\lambda }\left( \int_{\mathbb{R}%
^{N}}\left| u\right| ^{2_{\alpha }^{*}}dx\right) ^{1-\lambda }\leq C\frac{%
\left\| u\right\| _{A}^{p}}{A^{\frac{\left( 2N-2\right) \left( 1-\lambda
\right) }{2N-2-\alpha }}} 
\]
where $C:=S_{N}^{\lambda 2^{*}}\left( 2/\sigma _{N}\right) ^{\frac{2\alpha
\left( 1-\lambda \right) }{2N-2-\alpha }}$ only depends on $N,\alpha ,p$.
Recalling condition (\ref{growth p}), this implies 
\[
\left| \int_{\mathbb{R}^{N}}F\left( u\right) dx\right| \leq M_{2}\int_{\mathbb{R}%
^{N}}\left| u\right| ^{p}dx\leq M_{2}C\frac{\left\| u\right\| _{A}^{p}}{A^{%
\frac{\left( 2N-2\right) \left( 1-\lambda \right) }{2N-2-\alpha }}} 
\]
and therefore 
\[
I\left( u\right) \geq \frac{1}{2}\left\| u\right\| _{A}^{2}-a\left\|
u\right\| _{A}^{p} 
\]
where we set $a=M_{2}CA^{-\frac{\left( 2N-2\right) \left( 1-\lambda \right) 
}{2N-2-\alpha }}$ for brevity. Hence $I\left( tu\right) \geq \frac{1}{2}%
t^{2}\left\| u\right\| _{A}^{2}-at^{p}\left\| u\right\|
_{A}^{p}=:g_{u}\left( t\right) $ for every $t\geq 0$. The function $%
g_{u}:\left[ 0,+\infty \right) \rightarrow \mathbb{R}$ attains its maximum in $%
t_{u}:=(ap)^{-1/(p-2)}/\left\| u\right\| _{A}$ and, since 
\[
1-\lambda =\frac{p-2^{*}}{2_{\alpha }^{*}-2^{*}}=\frac{(p-2^{*})\left(
2N-2-\alpha \right) \left( N-2\right) }{4\left( \alpha -2\right) \left(
N-1\right) }, 
\]
one computes 
\begin{eqnarray*}
g_{u}\left( t_{u}\right) &=&\left( \frac{1}{ap}\right) ^{2/(p-2)}\left( 
\frac{1}{2}-\frac{1}{p}\right) =\frac{p-2}{2p^{p/(p-2)}}\left( \frac{1}{a}%
\right) ^{2/(p-2)}=\frac{p-2}{2p^{p/(p-2)}}\left( \frac{A^{\frac{\left(
2N-2\right) \left( 1-\lambda \right) }{2N-2-\alpha }}}{M_{2}C}\right)
^{2/(p-2)} \\
&=&\frac{p-2}{2p^{p/(p-2)}}\frac{A^{\frac{N-2}{\alpha -2}\frac{p-2^{*}}{p-2}}%
}{\left( M_{2}C\right) ^{\frac{2}{p-2}}}.
\end{eqnarray*}
Hence 
\[
\max_{t\geq 0}I\left( tu\right) \geq \max_{t\geq 0}g_{u}\left( t\right)
=g_{u}\left( t_{u}\right) =C_{0}A^{\frac{N-2}{\alpha -2}\frac{p-2^{*}}{p-2}} 
\]
with obvious definition of $C_{0}$. Since $u\in H_{\mathrm{r}}\setminus
\left\{ 0\right\} $ is arbitrary, we conclude 
\[
m_{A}=\inf_{u\in H_{\mathrm{r}}\setminus \left\{ 0\right\} ,\,u\geq
0}\max_{t\geq 0}I\left( tu\right) \geq \inf_{u\in H_{\mathrm{r}}\setminus
\left\{ 0\right\} }\max_{t\geq 0}I\left( tu\right) \geq C_{0}A^{\frac{N-2}{%
\alpha -2}\frac{p-2^{*}}{p-2}} 
\]
and the proof is complete.
\endproof

\begin{remark}
\label{RMK: p}If $p$ is as in Proposition \ref{PROP: mr>}, it is easy to
check that 
\[
\frac{N-2}{\alpha -2}\frac{p-2^{*}}{p-2}=\left\{ 
\begin{array}{lll}
\min \left\{ \frac{N-1}{\alpha },\frac{N-2}{2-\alpha }\frac{2^{*}-p_{1}}{%
p_{1}-2}\right\} \medskip  & \text{if } & 0<\alpha <2 \\ 
\min \left\{ \frac{N-1}{\alpha },\frac{N-2}{\alpha -2}\frac{p_{2}-2^{*}}{%
p_{2}-2}\right\}  & \text{if } & 2<\alpha <2N-2.
\end{array}
\right. 
\]
\end{remark}

\section{Estimate of $c_{A,K}$ \label{SEC: cs}}

Let $N\geq 3$, $2\leq K\leq N-2$ and $\alpha >0$, $\alpha \neq 2$. Let $f:%
\mathbb{R}\rightarrow \mathbb{R}$ be a continuous function satisfying $(\mathbf{f}%
_{p_{1},p_{2}})$, $(\mathbf{f}_{1})$, $(\mathbf{f}_{2})$. In this section we
define a suitable $\overline{u}_{K}\in H_{K}$ such that $I\left( \overline{u}%
_{K}\right) <0$ and estimate the corresponding mountain-pass level (\ref
{cs:=}).

In defining $\overline{u}_{K}$, we will use the following construction of
positive $H_{K}$ functions, which is inspired by \cite{Bad-Serra-mult}.
Denote by $\phi :D\rightarrow \mathbb{R}^{2}\setminus \left\{ 0\right\} $ the
change to polar coordinates in $\mathbb{R}^{2}\setminus \left\{ 0\right\} $,
namely $\phi \left( \rho ,\theta \right) =\left( \rho \cos \theta ,\rho \sin
\theta \right) $ for all $\left( \rho ,\theta \right) \in D:=\left(
0,+\infty \right) \times \left[ 0,2\pi \right) $. Define 
\[
E:=\left( \frac{1}{4},\frac{3}{4}\right) \times \left( \frac{\pi }{6},\frac{%
\pi }{3}\right) 
\]
and take any $\psi :\mathbb{R}^{2}\rightarrow \mathbb{R}$ such that $\psi \in
C_{c}^{\infty }\left( E\right) $ and $\psi >0$. For $0<\varepsilon <1$ and $%
\left( \rho ,\theta \right) \in \mathbb{R}^{2}$ define 
\[
\psi _{\varepsilon }\left( \rho ,\theta \right) :=\psi \left( \rho
^{1/\varepsilon },\frac{\theta }{\varepsilon }\right) , 
\]
in such a way that $\psi _{\varepsilon }\in C_{c}^{\infty }\left(
E_{\varepsilon }\right) $ where 
\[
E_{\varepsilon }:=\left\{ \left( \rho ,\theta \right) \in \mathbb{R}^{2}:\left(
\rho ^{1/\varepsilon },\frac{\theta }{\varepsilon }\right) \in E\right\}
=\left\{ \left( \rho ,\theta \right) \in \mathbb{R}^{2}:\left( \frac{1}{4}%
\right) ^{\varepsilon }<\rho <\left( \frac{3}{4}\right) ^{\varepsilon },%
\frac{\pi \varepsilon }{6}<\theta <\frac{\pi \varepsilon }{3}\right\} . 
\]
Finally define 
\[
v_{\varepsilon }\left( y,z\right) :=\psi _{\varepsilon }\left( \phi
^{-1}\left( \left| y\right| ,\left| z\right| \right) \right) \quad \text{for 
}x=\left( y,z\right) \in (\mathbb{R}^{K}\times \mathbb{R}^{N-K})\setminus \left\{
0\right\} ,\quad v_{\varepsilon }\left( 0\right) :=0. 
\]
Then $v_{\varepsilon }\in C_{c}^{\infty }\left( \Omega _{\varepsilon
}\right) \cap H_{K}$, where $\Omega _{\varepsilon }:=\left\{ \left(
y,z\right) \in \mathbb{R}^{K}\times \mathbb{R}^{N-K}:\left( \left| y\right|
,\left| z\right| \right) \in \phi \left( E_{\varepsilon }\right) \right\} $.

For future reference, we now compute the relevant integrals of $%
v_{\varepsilon }$. By means of spherical coordinates in $\mathbb{R}^{K}$ and $%
\mathbb{R}^{N-K}$ one has 
\begin{eqnarray*}
\int_{\mathbb{R}^{N}}\frac{v_{\varepsilon }^{2}}{\left| x\right| ^{\alpha }}dx
&=&\int_{\Omega _{\varepsilon }}\frac{\psi _{\varepsilon }\left( \phi
^{-1}\left( \left| y\right| ,\left| z\right| \right) \right) ^{2}}{\left|
x\right| ^{\alpha }}dx=\sigma _{K}\sigma _{N-K}\int_{\phi \left(
E_{\varepsilon }\right) }\frac{\psi _{\varepsilon }\left( \phi ^{-1}\left(
s,t\right) \right) ^{2}}{\left( s^{2}+t^{2}\right) ^{\alpha /2}}%
s^{K-1}t^{N-K-1}ds\,dt \\
&=&\sigma _{K}\sigma _{N-K}\int_{E_{\varepsilon }}\frac{\psi _{\varepsilon
}\left( \rho ,\theta \right) ^{2}}{\rho ^{\alpha -N+1}}H\left( \theta
\right) d\rho \,d\theta =\sigma _{K}\sigma _{N-K}\int_{E_{\varepsilon }}%
\frac{\psi \left( \rho ^{1/\varepsilon },\theta /\varepsilon \right) ^{2}}{%
\rho ^{\alpha -N+1}}H\left( \theta \right) d\rho \,d\theta
\end{eqnarray*}
where $H\left( \theta \right) :=\left( \cos \theta \,\right) ^{K-1}\left(
\sin \theta \right) ^{N-K-1}$, and by the change of variables 
\[
r=\rho ^{1/\varepsilon },\quad \varphi =\frac{\theta }{\varepsilon } 
\]
one obtains 
\begin{eqnarray}
\int_{\mathbb{R}^{N}}\frac{v_{\varepsilon }^{2}}{\left| x\right| ^{\alpha }}dx
&=&\sigma _{K}\sigma _{N-K}\int_{E}\frac{\psi \left( r,\varphi \right) ^{2}}{%
r^{\left( \alpha -N+1\right) \varepsilon }}H(\varepsilon \varphi
)\varepsilon ^{2}r^{\varepsilon -1}dr\,d\varphi  \nonumber \\
&=&\sigma _{K}\sigma _{N-K}\varepsilon ^{2}\int_{E}\frac{\psi \left(
r,\varphi \right) ^{2}}{r^{\left( \alpha -N\right) \varepsilon +1}}%
H(\varepsilon \varphi )dr\,d\varphi .  \label{comp1}
\end{eqnarray}
Similarly (recall that $F\left( 0\right) =0$) 
\begin{eqnarray}
\int_{\mathbb{R}^{N}}F\left( v_{\varepsilon }\right) dx &=&\int_{\Omega
_{\varepsilon }}F\left( \psi _{\varepsilon }\left( \phi ^{-1}\left( \left|
y\right| ,\left| z\right| \right) \right) \right) dx
\nonumber \\
&=&\sigma _{K}\sigma
_{N-K}\int_{\phi \left( E_{\varepsilon }\right) }F\left( \psi _{\varepsilon
}\left( \phi ^{-1}\left( s,t\right) \right) \right) s^{K-1}t^{N-K-1}ds\,dt 
\nonumber \\
&=&\sigma _{K}\sigma _{N-K}\int_{E_{\varepsilon }}F\left( \psi _{\varepsilon
}\left( \rho ,\theta \right) \right) \rho ^{N-1}H\left( \theta \right) d\rho
\,d\theta  \nonumber \\
&=&\sigma _{K}\sigma _{N-K}\varepsilon ^{2}\int_{E}F\left( \psi \left(
r,\varphi \right) \right) r^{N\varepsilon -1}H\left( \varepsilon \varphi
\right) dr\,d\varphi  \label{comp2}
\end{eqnarray}
and 
\begin{eqnarray}
\int_{\mathbb{R}^{N}}\left| \nabla v_{\varepsilon }\right| ^{2}dx
&=&\int_{\Omega _{\varepsilon }}\left| \nabla \psi _{\varepsilon }\left(
\phi ^{-1}\left( \left| y\right| ,\left| z\right| \right) \right) \cdot
J_{\phi ^{-1}}\left( \left| y\right| ,\left| z\right| \right) \right| ^{2}dx=
\nonumber \\
&=&\sigma _{K}\sigma _{N-K}\int_{\phi \left( E_{\varepsilon }\right) }\left|
\nabla \psi _{\varepsilon }\left( \phi ^{-1}\left( s,t\right) \right) \cdot
J_{\phi ^{-1}}\left( s,t\right) \right| ^{2}s^{K-1}t^{N-K-1}ds\,dt  \nonumber
\\
&=&\sigma _{K}\sigma _{N-K}\int_{E_{\varepsilon }}\left| \nabla \psi
_{\varepsilon }\left( \rho ,\theta \right) \cdot J_{\phi }^{-1}\left( \rho
,\theta \right) \right| ^{2}\rho ^{N-1}H\left( \theta \right) d\rho \,d\theta
\nonumber \\
&=&\sigma _{K}\sigma _{N-K}\int_{E_{\varepsilon }}\left( \frac{\partial \psi
_{\varepsilon }}{\partial \rho }\left( \rho ,\theta \right) ^{2}+\frac{1}{%
\rho ^{2}}\frac{\partial \psi _{\varepsilon }}{\partial \theta }\left( \rho
,\theta \right) ^{2}\right) \rho ^{N-1}H\left( \theta \right) d\rho \,d\theta
\nonumber \\
&=&\sigma _{K}\sigma _{N-K}\int_{E_{\varepsilon }}\frac{1}{\varepsilon ^{2}}%
\left( \rho ^{2/\varepsilon }\frac{\partial \psi }{\partial r}\left( \rho
^{1/\varepsilon },\frac{\theta }{\varepsilon }\right) ^{2}+\frac{\partial
\psi }{\partial \varphi }\left( \rho ^{1/\varepsilon },\frac{\theta }{%
\varepsilon }\right) ^{2}\right) \rho ^{N-3}H\left( \theta \right) d\rho
\,d\theta  \nonumber \\
&=&\sigma _{K}\sigma _{N-K}\int_{E}\left( \psi _{r}\left( r,\varphi \right)
^{2}+\frac{1}{r^{2}}\psi _{\varphi }\left( r,\varphi \right) ^{2}\right)
r^{(N-2)\varepsilon +1}H\left( \varepsilon \varphi \right) dr\,d\varphi .
\label{comp3}
\end{eqnarray}
where we denote $\psi _{r}=\frac{\partial \psi }{\partial r}$ and $\psi
_{\varphi }=\frac{\partial \psi }{\partial \varphi }$ for brevity.

\begin{lemma}
\label{LEM: lim+oo}The mapping $w_{A}:=v_{A^{-1/2}}\in H_{K}$, $A>1$, is
such that 
\[
\lim_{A\rightarrow +\infty }\frac{\left\| w_{A}\right\| _{A}^{2}}{\int_{\mathbb{%
R}^{N}}F\left( w_{A}\right) dx}=+\infty .
\]
\end{lemma}

\proof
According to the previous computations, for $\varepsilon =A^{-1/2}<1$ we
have 
\begin{eqnarray}
\frac{\left\| v_{\varepsilon }\right\| _{A}^{2}}{\int_{\mathbb{R}^{N}}F\left(
v_{\varepsilon }\right) dx} &=&\frac{\int_{\mathbb{R}^{N}}\left| \nabla
v_{\varepsilon }\right| ^{2}dx+A\int_{\mathbb{R}^{N}}\frac{v_{\varepsilon }^{2}%
}{\left| x\right| ^{\alpha }}dx}{\int_{\mathbb{R}^{N}}F\left( v_{\varepsilon
}\right) dx}  \nonumber \\
&=&\frac{\int_{E}\left( \left( \psi _{r}^{2}+\frac{1}{r^{2}}\psi _{\varphi
}^{2}\right) r^{(N-2)\varepsilon +1}+A\varepsilon ^{2}\psi ^{2}r^{\left(
N-\alpha \right) \varepsilon -1}\right) H(\varepsilon \varphi )dr\,d\varphi 
}{\varepsilon ^{2}\int_{E}F\left( \psi \right) r^{N\varepsilon -1}H\left(
\varepsilon \varphi \right) dr\,d\varphi }  \nonumber \\
&=&A\frac{\int_{E}\left( \left( \psi _{r}^{2}+\frac{1}{r^{2}}\psi _{\varphi
}^{2}\right) r^{(N-2)\varepsilon +1}+\psi ^{2}r^{\left( N-\alpha \right)
\varepsilon -1}\right) H(\varepsilon \varphi )dr\,d\varphi }{\int_{E}F\left(
\psi \right) r^{N\varepsilon -1}H\left( \varepsilon \varphi \right)
dr\,d\varphi }.  \label{R}
\end{eqnarray}
In the integration set $E$ one has $\varepsilon \pi /6<\varepsilon \varphi
<\varepsilon \pi /3$ and thus, for $\varepsilon >0$ small enough (i.e. $A>1$
large enough), we get that $\varepsilon \varphi /2<\sin \varepsilon \varphi
<\varepsilon \varphi $ and $1/2<\cos \varepsilon \varphi <1$. Hence there
exist two constants $\overline{C}_{1},\overline{C}_{2}>0$ such that 
\[
\overline{C}_{1}\varepsilon ^{N-K-1}<H\left( \varepsilon \varphi \right) <%
\overline{C}_{2}\varepsilon ^{N-K-1}. 
\]
Similarly, since $1/4<r<3/4$ in $E$, all the terms $r^{(N-2)\varepsilon +1}$%
, $r^{\left( N-\alpha \right) \varepsilon -1}$ and $r^{N\varepsilon -1}$ are
bounded and bounded away from zero by positive constants independent of $%
\varepsilon \in \left( 0,1\right) $ (i.e. of $A>1$), say $\overline{C}_{3}$
and $\overline{C}_{4}$ respectively. Inserting into (\ref{R}), this implies 
\begin{eqnarray*}
\frac{\left\| v_{\varepsilon }\right\| _{A}^{2}}{\int_{\mathbb{R}^{N}}F\left(
v_{\varepsilon }\right) dx} &\geq &A\frac{\overline{C}_{1}\int_{E}\left(
\left( \psi _{r}^{2}+\frac{1}{r^{2}}\psi _{\varphi }^{2}\right)
r^{(N-2)\varepsilon +1}+\psi ^{2}r^{\left( N-\alpha \right) \varepsilon
-1}\right) dr\,d\varphi }{\overline{C}_{2}\int_{E}Fr^{N\varepsilon
-1}dr\,d\varphi } \\
&\geq &A\frac{\overline{C}_{1}\overline{C}_{4}\int_{E}\left( \psi _{r}^{2}+%
\frac{1}{r^{2}}\psi _{\varphi }^{2}+\psi ^{2}\right) dr\,d\varphi }{%
\overline{C}_{2}\overline{C}_{3}\int_{E}F\left( \psi \right) dr\,d\varphi }.
\end{eqnarray*}
The last ratio is positive and independent of $A$, whence the claim follows.
\endproof

According to Lemma \ref{LEM: lim+oo}, we fix $A_{0}>1$ such that 
\begin{equation}
\frac{\left\| w_{A}\right\| _{A}^{2}}{\int_{\mathbb{R}^{N}}F\left( w_{A}\right)
dx}>1\quad \text{for every }A>A_{0}.  \label{def A0}
\end{equation}
We now distinguish the cases $0<\alpha <2$ and $\alpha >2$.

\begin{proposition}
\label{PROP: cs <2}Assume $(\mathbf{f}_{4})$ and $0<\alpha <2$. Let $A>A_{0}$
and define $\overline{u}_{K}\in H_{K}$ by setting 
\[
\overline{u}_{K}\left( x\right) :=w_{A}\left( \frac{x}{\lambda }\right)
\quad \text{with\quad }\lambda :=\frac{\left\| w_{A}\right\| _{A}^{2/\alpha }%
}{\left( \int_{\mathbb{R}^{N}}F\left( w_{A}\right) dx\right) ^{1/\alpha }}.
\]
Then $I\left( \overline{u}_{K}\right) <0$ and the corresponding
mountain-pass level (\ref{cs:=}) satisfies 
\[
c_{A,K}\leq C_{1}A^{\frac{K-1}{2}+N\left( \frac{1}{\alpha }-\frac{1}{2}%
\right) }
\]
where the constant $C_{1}>0$ does not depend on $A$.
\end{proposition}

\proof
Since $A>A_{0}$, one has $\lambda >1$. Then an obvious change of variables
yields 
\begin{eqnarray*}
I\left( \overline{u}_{K}\right) &=&\frac{\lambda ^{N-2}}{2}\int_{\mathbb{R}%
^{N}}\left| \nabla w_{A}\right| ^{2}dx+\frac{\lambda ^{N-\alpha }}{2}\int_{%
\mathbb{R}^{N}}\frac{A}{\left| x\right| ^{\alpha }}w_{A}^{2}dx-\lambda
^{N}\int_{\mathbb{R}^{N}}F\left( w_{A}\right) dx \\
&\leq &\frac{\lambda ^{N-\alpha }}{2}\left( \int_{\mathbb{R}^{N}}\left| \nabla
w_{A}\right| ^{2}dx+\int_{\mathbb{R}^{N}}\frac{A}{\left| x\right| ^{\alpha }}%
w_{A}^{2}dx\right) -\lambda ^{N}\int_{\mathbb{R}^{N}}F\left( w_{A}\right) dx \\
&=&\frac{\lambda ^{N}}{2}\left( \lambda ^{-\alpha }\left\| w_{A}\right\|
_{A}^{2}-2\int_{\mathbb{R}^{N}}F\left( w_{A}\right) dx\right) =-\frac{\lambda
^{N}}{2}\int_{\mathbb{R}^{N}}F\left( w_{A}\right) dx<0
\end{eqnarray*}
where the last inequality follows from assumption $(\mathbf{f}_{2})$, since $%
w_{A}>0$ almost everywhere. In order to estimate $c_{A,K}$, consider the
straight path $\gamma \left( t\right) :=t\overline{u}_{K}$, $t\in \left[
0,1\right] $. Clearly $c_{A,K}\leq \max_{t\in \left[ 0,1\right] }I\left(
\gamma \left( t\right) \right) $. Thanks to assumption $(\mathbf{f}_{4})$,
which implies $F\left( ts\right) \geq t^{\mu }F\left( s\right) $ for all $%
s>0 $ and $t\in \left[ 0,1\right] $, we have 
\[
I\left( \gamma \left( t\right) \right) =\frac{1}{2}t^{2}\left\| \overline{u}%
_{K}\right\| _{A}^{2}-\int_{\mathbb{R}^{N}}F\left( t\overline{u}_{K}\right)
dx\leq \frac{1}{2}t^{2}\left\| \overline{u}_{K}\right\| _{A}^{2}-t^{\mu
}\int_{\mathbb{R}^{N}}F\left( \overline{u}_{K}\right) dx=\frac{1}{2}%
t^{2}a-t^{\mu }b 
\]
where $a:=\left\| \overline{u}_{K}\right\| _{A}^{2}$ and $b:=\int_{\mathbb{R}%
^{N}}F\left( \overline{u}_{K}\right) dx$ for brevity. The function $g\left(
t\right) :=\frac{1}{2}t^{2}a-t^{\mu }b$ reaches its maximum in $t=\left( 
\frac{a}{b\mu }\right) ^{1/(\mu -2)}$, so that we get 
\[
I\left( \gamma \left( t\right) \right) \leq g\left( \left( \frac{a}{b\mu }%
\right) ^{1/(\mu -2)}\right) =a\left( \frac{a}{b\mu }\right) ^{2/(\mu
-2)}\left( \frac{1}{2}-\frac{1}{\mu }\right) . 
\]
Hence, setting $m:=\left( 1/\mu \right) ^{2/(\mu -2)}\left( 1/2-1/\mu
\right) $ for brevity and recalling that $\lambda >1$, we obtain 
\begin{eqnarray}
c_{A,K} &\leq &m\frac{\left\| \overline{u}_{K}\right\| _{A}^{2\mu /(\mu -2)}%
}{\left( \int_{\mathbb{R}^{N}}F\left( \overline{u}_{K}\right) dx\right)
^{2/(\mu -2)}}=m\frac{\left( \lambda ^{N-2}\int_{\mathbb{R}^{N}}\left| \nabla
w_{A}\right| ^{2}dx+\lambda ^{N-\alpha }\int_{\mathbb{R}^{N}}A\left| x\right|
^{-\alpha }w_{A}^{2}dx\right) ^{\mu /(\mu -2)}}{\left( \lambda ^{N}\int_{%
\mathbb{R}^{N}}F\left( w_{A}\right) dx\right) ^{2/(\mu -2)}}  \nonumber \\
&\leq &m\frac{\lambda ^{\mu \left( N-\alpha \right) /(\mu -2)}\left( \int_{%
\mathbb{R}^{N}}\left| \nabla w_{A}\right| ^{2}dx+\int_{\mathbb{R}^{N}}A\left|
x\right| ^{-\alpha }w_{A}^{2}dx\right) ^{\mu /(\mu -2)}}{\lambda ^{2N/(\mu
-2)}\left( \int_{\mathbb{R}^{N}}F\left( w_{A}\right) dx\right) ^{2/(\mu -2)}} 
\nonumber \\
&=&m\lambda ^{\frac{\mu \left( N-\alpha \right) -2N}{\mu -2}}\frac{\left\|
w_{A}\right\| _{A}^{2\mu /(\mu -2)}}{\left( \int_{\mathbb{R}^{N}}F\left(
w_{A}\right) dx\right) ^{2/(\mu -2)}}.  \label{PROP cs: c1}
\end{eqnarray}
Inserting the definition of $\lambda $ into (\ref{PROP cs: c1}), we get 
\[
c_{A,K}\leq m\frac{\left\| w_{A}\right\| _{A}^{\frac{2\mu }{\mu -2}+\frac{2}{%
\alpha }\frac{\mu \left( N-\alpha \right) -2N}{\mu -2}}}{\left( \int_{\mathbb{R}%
^{N}}F\left( w_{A}\right) dx\right) ^{\frac{2}{\mu -2}+\frac{1}{\alpha }%
\frac{\mu \left( N-\alpha \right) -2N}{\mu -2}}}=m\frac{\left\|
w_{A}\right\| _{A}^{\frac{2N}{\alpha }}}{\left( \int_{\mathbb{R}^{N}}F\left(
w_{A}\right) dx\right) ^{\frac{N-\alpha }{\alpha }}} 
\]
and therefore, using computations (\ref{comp1})-(\ref{comp3}) with $%
\varepsilon =A^{-1/2}$, we have 
\[
\overline{C}_{1}\varepsilon ^{N-K-1}<H\left( \varepsilon \varphi \right) <%
\overline{C}_{2}\varepsilon ^{N-K-1}. 
\]
\[
c_{A,K}\leq m\sigma _{K}\sigma _{N-K}\frac{\left( \int_{E}\left( \left( \psi
_{r}^{2}+\frac{1}{r^{2}}\psi _{\varphi }^{2}\right) r^{(N-2)\varepsilon
+1}+\psi ^{2}r^{\left( N-\alpha \right) \varepsilon -1}\right) H\left(
\varepsilon \varphi \right) dr\,d\varphi \right) ^{\frac{N}{\alpha }}}{%
\varepsilon ^{^{2\frac{N-\alpha }{\alpha }}}\left( \int_{E}F\left( \psi
\right) r^{N\varepsilon -1}H\left( \varepsilon \varphi \right) dr\,d\varphi
\right) ^{\frac{N-\alpha }{\alpha }}}. 
\]
As in the proof of Lemma \ref{LEM: lim+oo}, we take four constants $%
\overline{C}_{1},...,\overline{C}_{4}>0$ independent of $A$ such that for
every $\left( r,\varphi \right) \in E$ one has $\overline{C}_{1}\varepsilon
^{N-K-1}<H\left( \varepsilon \varphi \right) <\overline{C}_{2}\varepsilon
^{N-K-1}$ and the terms $r^{(N-2)\varepsilon +1}$, $r^{\left( N-\alpha
\right) \varepsilon -1}$ and $r^{N\varepsilon -1}$ are bounded and bounded
away from zero by $\overline{C}_{3}$ and $\overline{C}_{4}$ respectively.
Hence we conclude 
\begin{eqnarray*}
c_{A,K} &\leq &m\sigma _{K}\sigma _{N-K}\frac{\left( \overline{C}_{2}%
\overline{C}_{3}\int_{E}\left( \left( \psi _{r}^{2}+\frac{1}{r^{2}}\psi
_{\varphi }^{2}\right) +\psi ^{2}r\right) \varepsilon ^{N-K-1}dr\,d\varphi
\right) ^{\frac{N}{\alpha }}}{\varepsilon ^{^{2\frac{N-\alpha }{\alpha }%
}}\left( \overline{C}_{1}\overline{C}_{4}\int_{E}F\left( \psi \right)
\varepsilon ^{N-K-1}dr\,d\varphi \right) ^{\frac{N-\alpha }{\alpha }}} \\
&=&C\frac{\varepsilon ^{\frac{\left( N-K-1\right) N}{\alpha }}\left(
\int_{E}\left( \left( \psi _{r}^{2}+\frac{1}{r^{2}}\psi _{\varphi
}^{2}\right) +\psi ^{2}r\right) dr\,d\varphi \right) ^{\frac{N}{\alpha }}}{%
\varepsilon ^{^{\left( N-K+1\right) \frac{N-\alpha }{\alpha }}}\left(
\int_{E}F\left( \psi \right) dr\,d\varphi \right) ^{\frac{N-\alpha }{\alpha }%
}} \\
&=&CA^{\frac{K-1}{2}+N\left( \frac{1}{\alpha }-\frac{1}{2}\right) }\frac{%
\left( \int_{E}\left( \left( \psi _{r}^{2}+\frac{1}{r^{2}}\psi _{\varphi
}^{2}\right) +\psi ^{2}r\right) dr\,d\varphi \right) ^{\frac{N}{\alpha }}}{%
\left( \int_{E}F\left( \psi \right) dr\,d\varphi \right) ^{\frac{N-\alpha }{%
\alpha }}}
\end{eqnarray*}
with obvious definition of the constant $C$. As the last ratio does not
depend on $A$, the conclusion ensues.
\endproof

\begin{proposition}
\label{PROP: cs >2}Assume $(\mathbf{f}_{4})$ and $\alpha >2$. Let $A>A_{0}$
and define $\overline{u}\in H_{K}$ by setting 
\[
\overline{u}_{K}\left( x\right) :=w_{A}\left( \frac{x}{\lambda }\right)
\quad \text{with\quad }\lambda :=\frac{\left\| w_{A}\right\| _{A}}{\left(
\int_{\mathbb{R}^{N}}F\left( w_{A}\right) dx\right) ^{1/2}}.
\]
Then $I\left( \overline{u}_{K}\right) <0$ and the corresponding
mountain-pass level (\ref{cs:=}) satisfies 
\[
c_{A,K}\leq C_{2}A^{\frac{K-1}{2}}
\]
where the constant $C_{2}>0$ does not depend on $A$.
\end{proposition}

\proof
The proof is very similar to the one of Proposition \ref{PROP: cs <2}, so we
omit here some computational details. As $\alpha >2$, we have 
\begin{eqnarray*}
I\left( \overline{u}_{K}\right) &\leq &\frac{\lambda ^{N-2}}{2}\left( \int_{%
\mathbb{R}^{N}}\left| \nabla w_{A}\right| ^{2}dx+\int_{\mathbb{R}^{N}}\frac{A}{%
\left| x\right| ^{\alpha }}w_{A}^{2}dx\right) -\lambda ^{N}\int_{\mathbb{R}%
^{N}}F\left( w_{A}\right) dx \\
&=&-\frac{\lambda ^{N}}{2}\int_{\mathbb{R}^{N}}F\left( w_{A}\right) dx<0
\end{eqnarray*}
and 
\begin{eqnarray*}
\max_{t\in \left[ 0,1\right] }I\left( t\overline{u}_{K}\right) &\leq &m\frac{%
\left\| \overline{u}_{K}\right\| _{A}^{2\mu /(\mu -2)}}{\left( \int_{\mathbb{R}%
^{N}}F\left( \overline{u}_{K}\right) dx\right) ^{2/(\mu -2)}} \\
&=&m\frac{\left( \lambda ^{N-2}\int_{\mathbb{R}^{N}}\left| \nabla w_{A}\right|
^{2}dx+\lambda ^{N-\alpha }\int_{\mathbb{R}^{N}}A\left| x\right| ^{-\alpha
}w_{A}^{2}dx\right) ^{\mu /(\mu -2)}}{\left( \lambda ^{N}\int_{\mathbb{R}%
^{N}}F\left( w_{A}\right) dx\right) ^{2/(\mu -2)}} \\
&\leq &m\frac{\lambda ^{\mu \left( N-2\right) /(\mu -2)}\left( \int_{\mathbb{R}%
^{N}}\left| \nabla w_{A}\right| ^{2}dx+\int_{\mathbb{R}^{N}}A\left| x\right|
^{-\alpha }w_{A}^{2}dx\right) ^{\mu /(\mu -2)}}{\lambda ^{2N/(\mu -2)}\left(
\int_{\mathbb{R}^{N}}F\left( w_{A}\right) dx\right) ^{2/(\mu -2)}} \\
&=&m\lambda ^{\frac{\mu \left( N-2\right) -2N}{\mu -2}}\frac{\left\|
w_{A}\right\| _{A}^{2\mu /(\mu -2)}}{\left( \int_{\mathbb{R}^{N}}F\left(
w_{A}\right) dx\right) ^{2/(\mu -2)}}.
\end{eqnarray*}
Recalling the definition of $c_{A,K}$ and inserting the one of $\lambda $,
we get 
\[
c_{A,K}\leq m\frac{\left\| w_{A}\right\| _{A}^{\frac{2\mu }{\mu -2}+\frac{%
\mu \left( N-2\right) -2N}{\mu -2}}}{\left( \int_{\mathbb{R}^{N}}F\left(
w_{A}\right) dx\right) ^{\frac{2}{\mu -2}+\frac{1}{2}\frac{\mu \left(
N-2\right) -2N}{\mu -2}}}=m\frac{\left\| w_{A}\right\| _{A}^{N}}{\left(
\int_{\mathbb{R}^{N}}F\left( w_{A}\right) dx\right) ^{\frac{N-2}{2}}} 
\]
and therefore, using computations (\ref{comp1})-(\ref{comp3}) with $%
\varepsilon =A^{-1/2}$, we have 
\begin{eqnarray*}
c_{A,K} &\leq &m\sigma _{K}\sigma _{N-K}\frac{\left( \int_{E}\left( \left(
\psi _{r}^{2}+\frac{1}{r^{2}}\psi _{\varphi }^{2}\right) r^{(N-2)\varepsilon
+1}+\psi ^{2}r^{\left( N-\alpha \right) \varepsilon -1}\right) H\left(
\varepsilon \varphi \right) dr\,d\varphi \right) ^{\frac{N}{2}}}{\varepsilon
^{N-2}\left( \int_{E}F\left( \psi \right) r^{N\varepsilon -1}H\left(
\varepsilon \varphi \right) dr\,d\varphi \right) ^{\frac{N-2}{2}}} \\
&\leq &C\frac{\left( \int_{E}\left( \left( \psi _{r}^{2}+\frac{1}{r^{2}}\psi
_{\varphi }^{2}\right) +\psi ^{2}r\right) \varepsilon ^{N-K-1}dr\,d\varphi
\right) ^{\frac{N}{2}}}{\varepsilon ^{N-2}\left( \int_{E}F\left( \psi
\right) \varepsilon ^{N-K-1}dr\,d\varphi \right) ^{\frac{N-2}{2}}} \\
&=&CA^{\frac{K-1}{2}}\frac{\left( \int_{E}\left( \left( \psi _{r}^{2}+\frac{1%
}{r^{2}}\psi _{\varphi }^{2}\right) +\psi ^{2}r\right) dr\,d\varphi \right)
^{\frac{N}{2}}}{\left( \int_{E}F\left( \psi \right) dr\,d\varphi \right) ^{%
\frac{N-2}{2}}}
\end{eqnarray*}
where $C>0$ is a suitable constant independent of $A$. This concludes the
proof.
\endproof

\section{Proof of Theorem \ref{THM: main} \label{SEC: pf}}

This section is entirely devoted to the proof of Theorem \ref{THM: main}, so
we assume all the hypotheses of the theorem. The proof will be achieved
through some lemmas.

Let $K$ be any integer such that $2\leq K\leq N-2$. Assume $A>A_{K}$ (where $%
A_{K}$ is defined by (\ref{def A0})) and consider the mountain-pass level $%
c_{A,K}$ defined by (\ref{cs:=}), with $\overline{u}_{K}\in H_{K}$ given by
Lemma \ref{PROP: cs <2} or \ref{PROP: cs >2} 
according as $\alpha \in \left( 2/(N-1),2\right) $ or $\alpha \in \left(2,2N-2\right) $. 
We are going to show that $c_{A,K}$ is a critical level for the energy functional $I$ 
defined in (\ref{I(u)=}). To do this, we will make use of the sum space 
\[
L^{p_{1}}+L^{p_{2}}:=\left\{ u_{1}+u_{2}:u_{1}\in L^{p_{1}}\left( \mathbb{R}%
^{N}\right) ,\,u_{2}\in L^{p_{2}}\left( \mathbb{R}^{N}\right) \right\} . 
\]
We recall from \cite{BPR} that such a space can be characterized as the set
of measurable mappings $u:\mathbb{R}^{N}\rightarrow \mathbb{R}$ for which
there exists a measurable set $E\subseteq \mathbb{R}^{N}$ such that $u\in
L^{p_{1}}\left( E\right) \cap L^{p_{2}}(\mathbb{R}^{N}\setminus E)$ (%
\cite[Proposition 2.3]{BPR}). It is a Banach space with respect to the norm 
\[
\left\| u\right\| _{L^{p_{1}}+L^{p_{2}}}:=\inf_{u_{1}+u_{2}=u}\max \left\{
\left\| u_{1}\right\| _{L^{p_{1}}(\mathbb{R}^{N})},\left\| u_{2}\right\|
_{L^{p_{2}}(\mathbb{R}^{N})}\right\} 
\]
(\cite[Corollary 2.11]{BPR}) and the continuous embedding $L^{p}(\mathbb{R}%
^{N})\hookrightarrow L^{p_{1}}+L^{p_{2}}$ holds for all $p\in \left[
p_{1},p_{2}\right] $ (\cite[Proposition 2.17]{BPR}), in particular for $%
p=2^{*}$. Moreover, for every $u\in L^{p_{1}}+L^{p_{2}}$ and every $\varphi
\in L^{p_{1}^{\prime }}(\mathbb{R}^{N})\cap L^{p_{2}^{\prime }}(\mathbb{R}%
^{N})$ one has 
\begin{equation}
\int_{\mathbb{R}^{N}}\left| u\varphi \right| dx\leq \left\| u\right\|
_{L^{p_{1}}+L^{p_{2}}}\left( \left\| \varphi \right\| _{L^{p_{1}^{\prime }}(%
\mathbb{R}^{N})}+\left\| \varphi \right\| _{L^{p_{2}^{\prime }}(\mathbb{R}%
^{N})}\right)  \label{dual}
\end{equation}
where $p_{i}^{\prime }=p_{i}/(p_{i}-1)$ is the H\"{o}lder conjugate exponent
of $p_{i}$ (\cite[Lemma 2.9]{BPR}).

\begin{lemma}
\label{LEM: ck critical}$c_{A,K}$ is a critical level for the functional $%
I_{\mid H_{K}}$. 
\end{lemma}

\proof
Thanks to Lemma \ref{LEM: mp} (note that $I\left( \overline{u}_{K}\right) <0$
implies $\left\| \overline{u}_{K}\right\| _{A}>R$), the claim follows from
the Mountain Pass Theorem \cite{Ambr-Rab} if we show that $I_{\mid H_{K}}$
satisfies the Palais-Smale condition. Using the compact embeddings of \cite
{Azz-Pomp} and the results of \cite{BPR} about Nemytski\u{\i} operators on $%
L^{p_{1}}+L^{p_{2}}$, this is a standard proof but we still give some
details for the sake of completeness. Let $\left\{ u_{n}\right\} $ be a
sequence in $H_{K}$ such that $\left\{ I\left( u_{n}\right) \right\} $ is
bounded and $I^{\prime }\left( u_{n}\right) \rightarrow 0$ in the dual space
of $H_{K}$. Then, recalling (\ref{I(u)=}) and (\ref{I'(u)h=}), we have 
\[
\frac{1}{2}\left\| u_{n}\right\| _{A}^{2}-\int_{\mathbb{R}^{N}}F\left(
u_{n}\right) dx=O\left( 1\right) \quad \text{and}\quad \left\| u_{n}\right\|
_{A}^{2}-\int_{\mathbb{R}^{N}}f\left( u_{n}\right) u_{n}dx=o\left( 1\right)
\left\| u_{n}\right\| , 
\]
so that assumption $(\mathbf{f}_{1})$ implies 
\[
\frac{1}{2}\left\| u_{n}\right\| _{A}^{2}+O\left( 1\right) =\int_{\mathbb{R}%
^{N}}F\left( u_{n}\right) dx\leq \frac{1}{\theta }\int_{\mathbb{R}%
^{N}}f\left( u_{n}\right) u_{n}dx=\frac{1}{\theta }\left\| u_{n}\right\|
_{A}^{2}+o\left( 1\right) \left\| u_{n}\right\| . 
\]
This yields that $\left\{ \left\| u_{n}\right\| _{A}\right\} $ is bounded,
since $\theta >2$. On the other hand, thanks to the fact that $%
p_{1}<2^{*}<p_{2}$, the space $H_{K}$ is compactly embedded into $%
L^{p_{1}}+L^{p_{2}}$, since so is the subspace of $D^{1,2}(\mathbb{R}^{N})$
made up of the mappings with the same symmetries of $H_{K}$ (see 
\cite[Theorem A.1]{Azz-Pomp}). Hence there exists $u\in H_{K}$ such that, up
to a subsequence, we have $u_{n}\rightharpoonup u$ in $H_{K}$ and $%
u_{n}\rightarrow u$ in $L^{p_{1}}+L^{p_{2}}$. This implies that $\left\{
f\left( u_{n}\right) \right\} $ is bounded in both $L^{p_{1}^{\prime }}(%
\mathbb{R}^{N})$ and $L^{p_{2}^{\prime }}(\mathbb{R}^{N})$, since assumption 
$(\mathbf{f}_{p_{1},p_{2}})$ ensures that the operator $v\mapsto f\left(
v\right) $ is continuous from $L^{p_{1}}+L^{p_{2}}$ into $L^{p_{1}^{\prime
}}(\mathbb{R}^{N})\cap L^{p_{2}^{\prime }}(\mathbb{R}^{N})$ (see 
\cite[Corollary 3.7]{BPR}). Then by (\ref{dual}) we get 
\begin{eqnarray*}
\left| \int_{\mathbb{R}^{N}}f\left( u_{n}\right) \left( u_{n}-u\right)
dx\right| &\leq &\int_{\mathbb{R}^{N}}\left| f\left( u_{n}\right) \right|
\left| u_{n}-u\right| dx \\
&\leq &\left\| u_{n}-u\right\| _{L^{p_{1}}+L^{p_{2}}}\left( \left\| f\left(
u_{n}\right) \right\| _{L^{p_{1}^{\prime }}(\mathbb{R}^{N})}+\left\| f\left(
u_{n}\right) \right\| _{L^{p_{2}^{\prime }}(\mathbb{R}^{N})}\right) \\
&\leq &\left( \mathrm{const.}\right) \left\| u_{n}-u\right\|
_{L^{p_{1}}+L^{p_{2}}}=o\left( 1\right)
\end{eqnarray*}
and therefore 
\begin{eqnarray*}
\left\| u_{n}-u\right\| _{A}^{2} &=&\left( u_{n},u_{n}-u\right) _{A}-\left(
u,u_{n}-u\right) _{A} \\
&=&I^{\prime }\left( u_{n}\right) \left( u_{n}-u\right) +\int_{\mathbb{R}%
^{N}}f\left( u_{n}\right) \left( u_{n}-u\right) dx-\left( u,u_{n}-u\right)
_{A}=o\left( 1\right) ,
\end{eqnarray*}
where $\left( u,u_{n}-u\right) _{A}=o\left( 1\right) $ since $%
u_{n}\rightharpoonup u$ in $H_{K}$, and $I^{\prime }\left( u_{n}\right)
\left( u_{n}-u\right) =o\left( 1\right) $ because $I^{\prime }\left(
u_{n}\right) \rightarrow 0$ in the dual space of $H_{K}$ and $\left\{
u_{n}-u\right\} $ is bounded in $H_{K}$. This completes the proof.
\endproof

The next lemma clarifies why our separation of $c_{A,K}$ and $m_{A}$ needs
assumption (\ref{p1p2}) and the lower bound $\alpha >\frac{2}{N-1}$. Recall
the definition (\ref{n_alfa}) of $\nu =\nu _{N,\alpha ,p_{1},p_{2}}$.

\begin{lemma}
\label{LEM: conti}For every $\alpha \in \left( \frac{2}{N-1},2N-2\right) $, $%
\alpha \neq 2$, we have $\nu \geq 1$.
\end{lemma}

\proof
Assume $\frac{2}{N-1}<\alpha <2$. Since $\alpha >\frac{2}{N-1}$, we have 
\[
2\frac{N-1}{\alpha }-2N\left( \frac{1}{\alpha }-\frac{1}{2}\right) =N-\frac{2%
}{\alpha }>1. 
\]
On the other hand, by easy computations, condition 
\[
2\frac{N-2}{2-\alpha }\frac{2^{*}-p_{1}}{p_{1}-2}-2N\left( \frac{1}{\alpha }-%
\frac{1}{2}\right) >1 
\]
turns out to be equivalent to the first inequality of assumption (\ref{p1p2}%
). This proves that 
\[
2\min \left\{ \frac{N-1}{\alpha },\frac{N-2}{2-\alpha }\frac{2^{*}-p_{1}}{%
p_{1}-2}\right\} -2N\left( \frac{1}{\alpha }-\frac{1}{2}\right) >1, 
\]
which means 
\[
\left\lceil 2\min \left\{ \frac{N-1}{\alpha },\frac{N-2}{2-\alpha }\frac{%
2^{*}-p_{1}}{p_{1}-2}\right\} -2N\left( \frac{1}{\alpha }-\frac{1}{2}\right)
\right\rceil \geq 2 
\]
and thus $\nu \geq 1$. Similarly, if $2<\alpha <2N-2$ , we readily have $%
2(N-1)/\alpha >1$ and condition 
\[
2\frac{N-2}{\alpha -2}\frac{p_{2}-2^{*}}{p_{2}-2}>1 
\]
turns out to be equivalent to the second inequality of (\ref{p1p2}). This
proves again that $\nu \geq 1$.
\endproof

\proof[Proof of Theorem \ref{THM: main}]
On the one hand,
the restriction $I_{\mid H_{\mathrm{r}}}$ has a critical point $u_{\mathrm{r}%
}\neq 0$ thanks to the results of \cite{Su-Wang-Will-p}, since $(\mathbf{f}%
_{p_{1},p_{2}})$ ensures that one can find $p\in [p_{1},p_{2}]$ such that $%
\left| f\left( u\right) \right| \leq \left( \mathrm{const.}\right) u^{p-1}$
(cf. (\ref{growth p})) and (\ref{sww}) holds. On the other hand, according
to Lemma \ref{LEM: conti}, there are $\nu \geq 1$ integers $K$ (precisely $%
K=2,...,\nu +1$) such that 
\[
\frac{K-1}{2}+N\left( \frac{1}{\alpha }-\frac{1}{2}\right) <\min \left\{ 
\frac{N-1}{\alpha },\frac{N-2}{2-\alpha }\frac{2^{*}-p_{1}}{p_{1}-2}\right\}
\quad \text{if\quad }\frac{2}{N-1}<\alpha <2 
\]
and 
\[
\frac{K-1}{2}<\min \left\{ \frac{N-1}{\alpha },\frac{N-2}{\alpha -2}\frac{%
p_{2}-2^{*}}{p_{2}-2}\right\} \quad \text{if\quad }2<\alpha <2N-2. 
\]
Let $K$ be any of such integers. By Remark \ref{RMK: p} and Propositions \ref
{PROP: mr>}, \ref{PROP: cs <2} and \ref{PROP: cs >2}, there exists $%
A_{*}>A_{K}$ such that 
\begin{equation}
c_{A,K}<m_{A}\quad \text{for every }A>A_{*}\,.  \label{separation}
\end{equation}
Then, by Lemma \ref{LEM: ck critical}, there exists $u_{K}\in H_{K}$ such
that $I\left( u_{K}\right) =c_{A,K}$ and $I_{\mid H_{K}}^{\prime }\left(
u_{K}\right) =0$, where $u_{K}\neq 0$ since $c_{A,K}>0$ and $I\left(
0\right) =0$. Both $u_{\mathrm{r}}$ and $u_{K}$ are also critical points for
the functional $I:H_{\alpha }^{1}\rightarrow \mathbb{R}$, by the Palais'
Principle of Symmetric Criticality \cite{Palais}. Moreover, it easy to check
that they are nonnegative: test $I^{\prime }\left( u_{K}\right) $ with the
negative part $u_{K}^{-}\in H_{\alpha }^{1}$ of $u_{K}$ and use the fact
that $f\left( s\right) =0$ for $s<0$ to get $I^{\prime }\left( u_{K}\right)
u_{K}^{-}=-\left\| u_{K}^{-}\right\| _{A}^{2}=0$; the same for $u_{\mathrm{r}%
}$. Therefore $u_{\mathrm{r}}$ and $u_{K}$ are weak solutions to problem $%
\left( \mathcal{P}\right) $. Finally $u_{K}$ is not radial, because
otherwise Lemma \ref{LEM: mr} would imply $c_{A,K}=I\left( u_{K}\right) \geq
m_{A}$, which is false by (\ref{separation}). This also implies $%
u_{K_{1}}\neq u_{K_{2}}$ for $K_{1}\neq K_{2}$, thanks to Lemma \ref{LEM: Hk}.
\endproof

}

\end{document}